\documentclass[12pt]{article}
\usepackage[pctex32]{graphics}
\usepackage{amsthm,amsmath,amssymb,amscd,verbatim,epsfig}       % LaTex => PS => PDF
\usepackage{amssymb}
\usepackage{graphicx}
\date{}
\newcommand{\al}{\alpha}

\newcommand{\f}{\frac}
\newcommand{\ra}{\rightarrow}

\newcommand{\supp}{\supset}

\newcommand{\sq}{$\blacksquare$}

\begin{document}
\title{A\ simple\ proof\ of\ Sharkovsky's\ theorem\ rerevisited}
\author{Bau-Sen Du \\ %[.2cm]
Institute of Mathematics \\
Academia Sinica \\
Taipei 10617, Taiwan \\
dubs@math.sinica.edu.tw \\}
\maketitle
\begin{abstract}
Based on various strategies and a new general doubling operator, we present several simple directed-graph proofs of the celebrated Sharkovsky's cycle coexistence theorem.  A simple non-directed graph proof which is especially suitable for a calculus course right after the introduction of Intermediate Value Theorem is also given (in section 3).
\end{abstract}

\section{Introduction}
Throughout this note, $I$ is a compact interval, and $f : I \ra I$ is a continuous map.   For each integer $n \ge 1$, let $f^n$ be defined by: $f^1 = f$ and $f^n = f \circ f^{n-1}$ when $n \ge 2$.  For $x_0$ in $I$, we call $x_0$ a periodic point of $f$ with least period $m$ or a period-$m$ point of $f$ if $f^m(x_0) = x_0$ and $f^i(x_0) \ne x_0$ when $0 < i < m$.  If $f(x_0) = x_0$, then we call $x_0$ a fixed point of $f$.  

For discrete dynamical systems defined by iterated interval maps on $I$, one of the most unexpected results is Sharkovsky's cycle coexistence theorem which states as follows:   

\noindent
{\bf Theorem (Sharkovsky{\bf{\cite{mi, mi2, sh1, sh3}}})}  {\it Let the Sharkovsky's ordering of the natural numbers be defined (as suggested by Sharkovsky {\bf{\cite{sh3}}}) as follows: $$1 \prec 2 \prec 2^2 \prec 2^3 \prec \cdots \prec 2^n \prec \cdots \prec 9 \cdot 2^n \prec 7 \cdot 2^n \prec 5 \cdot 2^n \prec 3 \cdot 2^n \prec \cdots$$ $$\cdots \prec 9 \cdot 2^2 \prec 7 \cdot 2^2 \prec 5 \cdot 2^2 \prec 3 \cdot 2^2 \prec \cdots \prec 9 \cdot 2 \prec 7 \cdot2 \prec 5 \cdot 2 \prec 3 \cdot 2 \prec \cdots \prec 9 \prec 7 \prec 5 \prec 3.$$   Then the following three statements hold:
\begin{itemize}
\item[\rm{(1)}] 
If $f$ has a period-$m$ point and if $n \prec m$, then $f$ also has a period-$n$ point.

\item[\rm{(2)}]
For each positive integer $n$ there exists a continuous map from $I$ into itself that has a period-$n$ point but has no period-$m$ point for any $m$ with $n \prec m$.

\item[\rm{(3)}]
There exists a continuous map from $I$ into itself that has a period-$2^i$ point for $i = 0, 1, 2, \ldots$ but has no periodic point of any other period.
\end{itemize}}

It is clear that (1) is equivalent to the following three statements (cf. {\bf{\cite{du2, du3, str}}}): 
\begin{itemize}
\item[(a)] if $f$ has a period-$m$ point with $m = 3$ or $4$, then $f$ has a period-2 point and a fixed point; 

\item[(b)] if $f$ has a period-$m$ point with $m \ge 3$ and odd, then $f$ has a period-$(m+2)$ point; and 

\item[(c)] if $f$ has a period-$m$ point with $m \ge 3$ and odd, then $f$ has a period-$6$ point and a period-$(2m)$ point.  
\end{itemize}

\noindent
Note that in (a) we only need the special cases when $m = 3,4$ (instead of all $m \ge 3$) which can be easily proved.  Also, in (b) and (c), we don't require the existence of periodic points of all periods $\ge m+1$ and all {\it even} periods.  Only the existence of period-$(m+2)$, period-6 and period-$(2m)$ points suffices.  See section 11 for details.  However, the proofs of these special cases turn out to be no easier than those of the general ones.  They only serve to explain why Sharkovsky's ordering is defined as it is.  

In the past 35 years, there have been a number of papers dealing with Sharkovsky's theorem (see references), including the three papers {\bf{\cite{du1, du2, du3}}} by the author and the "standard proof" developed in {\bf{\cite{bc, bl, bu, ho, str}}} and improved in {\bf{\cite{al}}} which, for the odd period cases, shows the existence of \v Stefan cycles first and then draws conclusions on (1) of Sharkovsky's theorem from the directed graphs of such cycles.  On the other hand, in {\bf{\cite{bh}}}, Burns and Hasselblatt "select a salient sequence of orbit points and prove that this sequence `spirals out` in essentially the same way as the \v Stefan cycles considered in the standard proof".  Their proof and the standard's all use cycles of compact intervals with endpoints belonging to one and the same periodic orbit.  By allowing these endpoints to be periodic or preperiodic points in different orbits, we have more cycles at our disposal and so can reach the goal more easily.  Based on this tactic, we present in {\bf{\cite{du2}}} a simple proof of (1) of Sharkovsky's theorem by going around the \v Stefan cycles.  In this note, we shall present several strategies on how to prove (a), (b) and (c).  We even confront the \v Stefan cycles with some quite straightforward arguments which achieve the same goal as the standard proof and yet surprisingly are just as simple as that in {\bf{\cite{du2}}}.  

To make this note self-contained, we include some well-known preliminary results in section 2.  In section 3, we give a non-directed graph proof of (a), (b) and (c) which uses the Intermediate Value Theorem in a very straightforward way.  In sections 4, 5 and 6, we concentrate our attention on the point $\min P$ and/or the point $\max P$ in any given period-$m$ orbit $P$.  In section 7, we examine how the iterates of a given period-$m$ point "jump" around a fixed point of $f$.  In section 9, we investigate how the points in a given period-$m$ orbit which lie on either side of a fixed point of $f$ are mapped to the other side by $f$.  In these two sections, we even prove, when $m \ge 5$ is odd, the existence of \v Stefan cycles of least period $m$ without assuming the non-existence of smaller odd periods other than fixed points.  In sections 8 and 10, we use the strategy in each previous section to treat the special case when $f$ has a periodic point of odd period $m \ge 5$ but no periodic points of smaller odd periods other than fixed points.  The proof in each section is independent of the other.  Finally, in section 11, we present two proofs of Part (1) of Sharkovsky's theorem.  In section 12, we introduce a new general doubling operator which, together with the classical one as described in {\bf\cite{al, st}}, is used in various combinations to construct new examples for (2) and (3).  

\section{Preliminary results}
To make this paper self-contained, we include the following well-known results.

\noindent
{\bf Lemma 1.}
If $f^n(x_0) = x_0$, then the least period of $x_0$ with respect to $f$ divides $n$.

\noindent
{\it Proof.}
Let $m$ denote the least period of $x_0$ with respect to $f$ and write $n = km + r$ with $0 \le r < m$.  Then $x_0 = f^n(x_0) = f^{km+r}(x_0) = f^r(f^{km}(x_0)) = f^r(x_0)$.  Since $m$ is the smallest positive integer such that $f^m(x_0) = x_0$, we must have $r = 0$.  Therefore, $m$ divides $n$.  
\hfill\sq

\noindent
{\bf Lemma 2.}
If $J$ is a closed subinterval of $I$ and $f(J) \supset J$, then $f$ has a fixed point in $J$.

\noindent
{\it Proof.}
Write $J = [a, b]$.  Since $f(J) \supset J \supset \{ a, b \}$, there exist points $p$ and $q$ in $[a, b]$ such that $f(p) = a$ and $f(q) = b$.  Let $g : I \longrightarrow R$ be a continuous map defined by $g(x) = f(x) - x$.  Then $g(p) = f(p) - p = a - p \le 0$ and $g(q) = f(q) - q = b - q \ge 0$.  By Intermediate Value Theorem, there is a point $z$ between $p$ and $q$ such that $f(z) - z = g(z) = 0$.  So, $z$ is a fixed point of $f$ in $J$.  
\hfill\sq

\noindent
{\bf Lemma 3.}
{\it Let $k, m, n$, and $s$ be positive integers.  Then the following statements hold:
\begin{itemize}
\item[\rm{(1)}]
If $x_0$ is a periodic point of $f$ with least period $m$, then it is a periodic point of $f^n$ with least period $m/(m, n)$, where $(m, n)$ is the greatest common divisor of $m$ and $n$. 

\item[\rm{(2)}]
If $x_0$ is a periodic point of $f^n$ with least period $k$, then it is a periodic point of $f$ with least period $kn/s$, where $s$ divides $n$ and is relatively prime to $k$. 
\end{itemize}}

\noindent
{\it Proof.} (1) Let $t$ denote the least period of $x_0$ under $f^n$.  Then $m$ divides $nt$ since $x_0 = (f^n)^t(x_0) = f^{nt}(x_0)$.  Consequently, $\f m{(m, n)}$ divides $\f n{(m, n)} \cdot t$.  Since $\f m{(m, n)}$ and $\f n{(m, n)}$ are coprime, $\f m{(m, n)}$ divides $t$.  On the other hand, $(f^n)^{(m/(m, n))}(x_0) = (f^m)^{(n/(m, n))}(x_0) = x_0$.  So, $t$ divides $\f m{(m, n)}$.  This shows that $t = \f m{(m, n)}$.

(2) Since $x_0 = (f^n)^k(x_0) = f^{kn}(x_0)$, the least peirod of $x_0$ under $f$ is $\f {kn}s$ for some positive integer $s$.  By (1), $(\f {kn}s)/((\f {kn}s), n) = k$.  So, $\f ns = ((\f ns)k, n)$ (which is an integer) $= ((\f ns)k, (\f ns)s) = (\f ns)(k, s)$.  This shows that $s$ divides $n$ and $(s, k) = 1$.
\hfill\sq

\noindent
{\bf Lemma 4.}
Let $J$ and $L$ be closed subintervals of $I$ with $f(J) \supset L$.  Then there exists a closed subinterval $K$ of $J$ such that $f(K) = L$.

\noindent
{\it Proof.}
Let $L = [a, b]$.  Then since $\{ a, b \} \subset L \subset f(J)$, there are two points $p$ and $q$ in $J$ such that $f(p) = a$ and $f(q) = b$.  If $p < q$, let $c = \max \{ p \le x \le q : f(x) = a \}$ and let $d = \min \{ c \le x \le q : f(x) = b \}$.  If $p > q$, let $c = \max \{ q \le x \le p \mid f(x) = b \}$ and let $d = \min \{ c \le x \le p : f(x) = a \}$.  In either case, let $K = [c, d]$.  Then $f(K) = L$.
\hfill\sq

If there are closed subintervals $J_0$, $J_1$, $\cdots$, $J_{n-1}, J_n$ of $I$ with $J_n = J_0$ such that $f(J_i) \supp J_{i+1}$ for $i = 0, 1, \cdots, n-1$, then we say that $J_0J_1 \cdots J_{n-1}J_0$ is a {\it{cycle of length}} $n$.  We need  the following result which is useful for showing the existence of periodic points of some periods.  

\noindent
{\bf Lemma 5.}
{\it If $J_0J_1J_2 \cdots J_{n-1}J_0$ is a cycle of length $n$, then there exists a periodic point $y$ of $f$ such that $f^i(y)$ belongs to $J_i$ for $i = 0, 1, \cdots, n-1$ and $f^n(y) = y$.}

\noindent
{\it Proof.}
Let $Q_n = J_0$.  Since $f(J_{n-1}) \supset J_0 = Q_n$, there is, by Lemma 4, a closed subinterval $Q_{n-1}$ of $J_{n-1}$ such that $f(Q_{n-1}) = Q_n = J_0$.  Continuing this process one by one, we obtain, for each $0 \le i \le n-1$, a closed subinterval $Q_i$ of $J_i$ such that $f(Q_i) = Q_{i+1}$.  Consequently, $f^i(Q_0) = Q_i$ for all $0 \le i \le n$.  In particular, $f^n(Q_0) = Q_n = J_0 \supset Q_0$.  By Lemma 2, there is a point $y$ in $Q_0 \subset J_0$ such that $f^n(y) = y$.  Since $y \in Q_0$, we also obtain that $f^i(y) \in f^i(Q_0) = Q_i \subset J_i$ for all $0 \le i \le  n-1$. 
\hfill\sq

\noindent
{\bf Remark.}
The point $y$ obtained in Lemma 5 need not have least period $n$ in general.  However, by choosing appropriate cycles of length $n$, we can still get periodic points of least period $n$.  

\section{A non-directed graph proof of (a), (b) and (c)}
The proof we present here is a slight improvement of the one in {\bf{\cite{du3}}}.  It is more direct and assumes no knowledge of dynamical systems theory whatsoever.      

Let $P$ be a period-$m$ orbit of $f$ with $m \ge 2$ and let $b = f^{m-1}(\min P)$.  Then $f(b) = \min P < b$.  If $f(x) < b$ on  $[\min P, b]$, then, $(\min P \le) \,\, f^i(\min P) < b$ for all $i \ge 1$, contradicting the fact that $f^{m-1}(\min P) = b$.  So, there is a point $a$ in $[\min P, b]$ such that $f(a) \ge b$.  Thus, $f$ has a fixed point $z$ in $[a, b]$.  Now suppose $m \ge 3$ and let $v$ be a point in $[a, z]$ such that $f(v) = b$.  Since $f^2(\min P) > \min P$ and $f^2(v) = \min P < v$, the point $y = \max \{ \min P \le x \le v : f^2(x) = x \}$ exists.  Furthermore, $f(x) > z$ on $[y, v]$ and $f^2(x) < x$ on $(y, v]$.  Therefore, $y$ is a period-2 point of $f$.  (a) is proved.

For the proofs of (b) and (c), we assume that $m \ge 3$ is odd and note that $f(x) > z > x > f^2(x)$ on $(y, v]$.  Since $f^{m+2}(y) = f(y) > y$ and $f^{m+2}(v) = f^m(\min P) = \min P < v$, the point $p_{m+2} = \min \{ y \le x \le v : f^{m+2}(x) = x \}$ exists.  Let $k$ denote the least period of $p_{m+2}$ with respect to $f$.  Then $k > 1$ and, by Lemma 1, $k$ divides $m+2$ and so is odd.  If $k < m+2$, then since $f^{k+2}(y) = f(y) > y$ and $f^{k+2}(p_{m+2}) = (f^2)(f^k(p_{m+2})) = f^2(p_{m+2}) < p_{m+2}$, there is a point $w_{k+2}$ in $(y, p_{m+2})$ such that $f^{k+2}(w_{k+2}) = w_{k+2}$.  Inductively, there exist points $$y < \cdots < w_{m+2} < w_m < w_{m-2} < \cdots < w_{k+4} < w_{k+2} < p_{m+2} < v$$ such that $f^{k+2i}(w_{k+2i}) = w_{k+2i}$ for all $i \ge 1$.  In particular, $f^{m+2}(w_{m+2}) = w_{m+2}$ and $y < w_{m+2} < p_{m+2}$, contradicting the fact that $p_{m+2}$ is the {\it smallest} point in $(y, v)$ which satisfies $f^{m+2}(x) = x$.  Therefore, $k = m+2$.  This establishes (b).  

We now prove (c).  Let $z_0 = \min \{ v \le x \le z : f^2(x) = x \}$.  Then $f^2(x) < x < z_0$ and $f(x) > z$ on $(v, z_0)$ and so also on $(y, z_0)$.  If $f^2(x) < z_0$ whenever $\min P \le x < z_0$, then we have $\min P \le f^{2i}(\min P) < z_0$ for all $i \ge 1$ which contradicts the fact that $(f^2)^{(m-1)/2}(\min P) = b > z_0$.  Since $f^2(x) < x < z_0$ on $(y, z_0)$, the point $$d = \max \{ \min P \le x \le y : f^2(x) = z_0 \}$$ exists and $f(x) > z \ge z_0 > f^2(x)$ on $(d, y)$.  Therefore, $f(x) > z \ge z_0 > f^2(x)$ on $(d, z_0)$.  Let $$u_1 = \min \{ d \le x \le v : f^2(x) = d \}.$$  Then $d < f^2(x) < z_0$ on $(d, u_1)$.  Let $c_1$ be any point in $(d, u_1)$ such that $f^2(c_1) = c_1$.  The point $$u_2 = \min \{ d \le x \le c_1 : (f^2)^2(x) = d \}$$exists since $f^2([d, c_1]) \supset [c_1, z_0] \supset \{ u_1 \}$.  Then $d < (f^2)^2(x) < z_0$ on $(d, u_2)$.  Let $c_2$ be any point in $(d, u_2)$ such that $(f^2)^2(c_2) = c_2$.  Let $$u_3 = \min \{ d \le x \le c_2 : (f^2)^3(x) = d \}.$$  Then $d < (f^2)^3(x) < z_0$ on $(d, u_3)$.  Let $c_3$ be any point in $(d, u_3)$ such that $(f^2)^3(c_3) = c_3$.  Proceeding in this manner indefinitely, we obtain points $$d < \cdots < c_n < u_n < \cdots < c_2 < u_2 < c_1 < u_1 < z_0$$ such that $d < (f^2)^n(x) < z_0$ on $(d, u_n)$ and $(f^2)^n(c_n) = c_n$.  Since $f(x) > z \ge z_0$ on $(d, z_0)$, we have $$f^i(c_n) < z_0 < f^j(c_n) \,\,\, \text{for all even} \,\,\, i \,\,\, \text{and all odd} \,\,\, j \,\,\, \text{in} \,\,\, [0, 2n].$$  Therefore, each $c_n$ is a period-$(2n)$ point of $f$.  This proves (c).

\noindent
{\bf Remarks.}  (1) It is well-known (see the first paragraph of section 6 for a proof) that if there exist a fixed point $\hat z$ of $h \in C^0(I, I)$, a point $c$ in $I$ and an integer $n \ge 2$ such that $h(c) < c < \hat z \le h^n(c)$ then $h$ has periodic points of all periods.  In the above proof of (c), we have $f^2(v) < v < z < b = (f^2)^{(m+1)/2}(v)$.  So, $f^2$ has periodic points of all periods.  However, we cannot conclude the existence of period-$(2j)$ points for $f$ for any odd $j \ge 3$ from this fact yet.  We need to do a little more work to ensure that as we did above.   The following is another approach:  Let $m, v, z, b$ be defined as in the above proof.  Let $g \in C^0(I, I)$ be defined by $g(x) = \max \{ f(x), z \}$ if $x \le z$ and $g(x) = \min \{ f(x), z \}$ if $x \ge z$.  Then $g([\min I, z]) \subset [z, \max I]$ and $g([z, \max I]) \subset [\min I, z]$.  So, $g$ has {\it no} periodic points of any {\it odd} periods $\ge 3$.  Since $m \ge 3$ is odd, for some $1 \le i \le m-1$, both $f^i(b)$ and $f^{i+1}(b)$ lie on the same side of $z$.  Let $k$ be the {\it smallest} such $i$.  Then the iterates $b, f(b), f^2(b), \cdots, f^k(b)$ are jumping alternately around $z$ and since $f(v) = b$, so are the iterates $v, f(v), f^2(v), \cdots, f^{k+1}(v)$.  Consequently, $g^i(v) = f^i(v)$ for all $0 \le i \le k+1$.  If $k$ is odd, then $f^k(b) < z$ and $f^{k+1}(b) < z$ and so $g^{k+1}(b) = z$ and $g^2(v) = f^2(v) < v < z = g(z) = g^{k+2}(b) = g^{k+3}(v)$.  If $k \ge 2$ is even, then $f^k(b) > z$ and $f^{k+1}(b) > z$ and so $g^{k+1}(b) = z$ and $g^2(v) = f^2(v) < v < z = g^{k+2}(v)$.  In either case, we have $g^2(v) < v < z = (g^2)^n(v)$ for some $n \ge 2$.  As noted above, this implies that $g^2$ has periodic points of all periods.  In particular, $g^2$ has period-$j$ points for all odd $j \ge 3$.  So, $g$ has either period-$j$ points or period-$(2j)$ points for any odd $j \ge 3$.  Since $g$ has no periodic points of any odd periods $\ge 3$, $g$ has period-$(2j)$ points which are also period-$(2j)$ points of $f$ for all odd $j \ge 3$.  This establishes (c).

If for each odd integer $\ell \ge 3$, $f$ has only finitely many period-$\ell$ points, then we have the following different approach: Let $P$ be a period-$m$ orbit of $f$ such that $(\min P, \max P)$ contains no period-$m$ orbits of $f$.  Then by (b), $[\min P, \max P]$ contains a period-$(m+2)$ orbit $Q$ of $f$.  Let $h$ be the continuous map from $I$ into itself defined by $h(x) = \min Q$ if $f(x) \le \min Q$; $h(x) = \max Q$ if $f(x) \ge \max Q$; and $h(x) = f(x)$ elsewhere.  Then $h$ has no period-$m$ points and by (b) $h$ has no period-$j$ points for any odd $3 \le j \le m-2$.  Since $h$ has the period-$(m+2)$ orbit $Q$, there exist a fixed point $z$ of $h$ and a point $v$ such that $h(v) \in Q$ and  $h^2(v) < v < z < h(v) = (h^2)^{(m+3)/2}(v)$.  Consequently, $h^2$ has periodic points of all periods.  Since $h$ has no period-$j$ points for any odd $3 \le j \le m$, $h$ has period-$(2j)$ points for all odds $3 \le j \le m$ which are also periodic points of $f$ with the same periods.  In particular, $f$ has period-6 and period-$(2m)$ points and so (c) is proved.  

(2) The arguments in the above proofs of (a) and (b) can be used to give a simpler proof of the main result of Block in {\bf{\cite{bl0}}} on the stability of periodic orbits in Sharkovsky's theorem without resorting to the \v Stefan cycles.  Indeed, assume that $f$ has a period-$2^n$ point.  Let $F = f^{2^{n-2}}$.  Then $F$ has a period-4 orbit $Q$.  Arguing as in the proof of (a), there exists a point $v$ such that $\min Q = F^2(v) < v < F(v)$.  Since $F^2(\min Q) > \min Q$, there is an open neighborhood $U$ of $f$ in $C^0(I, I)$ such that, for each $g$ in $U$, the map $G = g^{2^{n-2}}$ satisfies $G^2(v) < v < G(v)$ and $G^2(\min Q) > \min Q$.  Thus, the point $y = \max \{ \min Q \le x \le v : G^2(x) = x \}$ is a period-2 point of $G$.  Consequently, $y$ is a period-$2^{n-1}$ point of $g$.  On the other hand, assume that $f$ has a period-$(2^i \cdot m)$ point with $m \ge 3$ and odd and $i \ge 0$.  Let $F = f^{2^i}$.  Then $F$ has period-$m$ points.  Arguing as in the proof of (a), there exist a period-2 point $y$ of $F$ and a point $v$ with $y < v$ such that $F^{m+2}(y) > y$ and $F^{m+2}(v) = F^2(v) < v$ and $F$ has no fixed points in $[y, v]$.  So, there is an open neighborhood $W$ of $f$ in $C^0(I, I)$ such that, for each $h$ in $W$, the map $H = h^{2^i}$ satisfies that $H^{m+2}(y) > y$, $H^{m+2}(v) < v$, and $H$ has no fixed points in $[y, v]$.  By Sharkovsky's theorem, $h$ has period-$(2^i \cdot (m+2))$ points.

From now on, let $P$ be a period-$m$ orbit of $f$ with $m \ge 3$.  We shall present seven different directed graph (digraph for short) proofs of (a), (b) and (c).  The first three proofs depend on the choice of the particular point $\min P$ and/or $\max P$.  These particular choices simplify the proofs of (a), (b) or even (c) sometimes.  In some cases, we shall need the following easy fact that if $f$ has a period-3 point then $f$ has a fixed point $z$ and, for each $n \ge 2$, has a period-$n$ point $p_n$ whose iterates $p_n, f(p_n), f^2(p_n), \cdots, f^{n-1}(p_n)$ are "spiralling out" alternately around the fixed point $z$ of $f$.  In particular, $f$ has periodic points of all periods.  We shall also need the following easy result to show the existence of periodic points of all {\it even} periods.  

\noindent
{\bf Lemma 6.}  {\it Assume that there exist a point $d$ and a fixed point $z$ of $f$ such that $$f^3(d) \le z \,\,\, \text{and} \,\,\, f^2(d) < d < z < f(d), \quad \text{\rm or} \quad f(d) < z < d < f^2(d) \,\,\, \text{and} \,\,\, z \le f^3(d).$$ Then $f$ has periodic points of all even periods.  Furthermore, these periodic points can be chosen between $z$ and $f(d)$ so that their iterates which lie between $f(d)$ and $f^2(d)$ "spiral out" alternately around the fixed point $z$.}

\noindent
{\it Proof.}
We may assume that $f^3(d) \le z$ and $f^2(d) < d < z < f(d)$.  Let $x_{-i}, i \ge 0$ be points such that $f^2(d) < d = x_0 < x_{-2} < x_{-4} < \cdots < z < \cdots < x_{-5} < x_{-3} < x_{-1} < f(d)$ and $f(x_{-j}) = x_{-j+1}$ for all $j \ge 1$.  By considering the cycles $[x_{-1}, f(d)][f^2(d), d][x_{-1}, f(d)]$ and, $$[x_{-2n+1}, x_{-2n+3}] [x_{-2n+4}, x_{-2n+2}] \cdots [x_{-3}, x_{-1}] [d, x_{-2}] [x_{-1}, f(d)] [f^2(d), d] [x_{-2n+1}, x_{-2n+3}], n \ge 2$$ of length $2n$, we obtain period-$(2n)$ points between $z$ and $f(d)$ whose iterates (which lie between $f(d)$ and $f^2(d)$) "spiral out" alternately around the fixed point $z$.

\section{The first digraph proof of (a), (b) and (c)}
Let $b = f^{m-1}(\min P)$.  Then $f(b) = \min P$.  If $f(x) < b$ for all $\min P \le x \le b$, then $\min P \le f^i(\min P) < b$ for all $i \ge 1$.  This contradicts the fact that $b = f^{m-1}(\min P)$.  Consequently, since $f(b) = \min P$, we have $f([\min P, b]) \supset [\min P, b]$.  Let $v$ be a point in $(\min P, b)$ such that $f(v) = b$ and let $z$ be a fixed point of $f$ in $(v, b)$.    

If for some point $x$ in $[\min P, v]$, $f(x) \le z$, then we consider the cycle $[x, v][z, b][x, v]$.  If for some $x$ in $[\min P, v]$, $f(x) > z$ and $v < f^2(x)$, then we consider the cycle $[x, v][b : f(x)][x, v]$, where $[c : d] = [c, d]$ if $c < d$ and $[c : d] = [d, c]$ if $c > d$.  If $f^2(x) \le v < z < f(x)$ for all $x$ in $[\min P, v]$, let $p = \max(P \cap [\min P, v]) \,\, (\le v)$.  Then $f^2(P \cap [\min P, p]) \subset P \cap [\min P, p]$.  Since $f^2$ is one-to-one on $P$, this implies that $f^2(P \cap [\min P, p]) = P \cap [\min P, p]$ (and so $m$ must be even).  Consequently, $f^2([\min P, p]) \supset [\min P, p]$.  In this case, we consider the cycle $[\min P, p]f([\min P, p])[\min P, p]$.  In either case, we obtain a period-2 point $y$ of $f$ such that $y < v < f(y)$ (this fact will be used below).  This proves (a).  

We now prove (b).  Let $m \ge 3$ be odd.  Let $k$ be the {\it largest even} integer in $[0, m-1]$ such that $f^k(\min P) < y$.  Then $0 \le k \le m-2$ since $f^{m-1}(\min P) = b > y$.  For $0 \le i \le m-2$, let $J_i = [f^i(\min P) : f^i(y)]$.  If there is a {\it smallest odd} integer $\ell$ with $k < \ell \le m-2$ such that $f^\ell(\min P) < y$, then $y < f^i(\min P)$ for all $k < i < \ell$.  We consider the cycle $$J_kJ_{k+1}J_{k+2} \cdots J_{\ell-1} [v, f(y)]([v, b])^{n+k-\ell-1}J_k$$ of length $n$ for each $n \ge m+1$, where $([v, b])^{n+k-\ell-1}$ denotes the $n+k-\ell-1$ copies of $[v, b]$.  If $y < f^j(\min P)$ for all odd integers $k < j \le m-2$, we consider the cycle $$J_kJ_{k+1}J_{k+2} \cdots J_{m-2}([v, b])^{n-m+k+1}J_k$$ of length $n$ for each $n \ge m+1$.  In either case, we obtain a period-$n$ point $p_n$ in $[\min P, y]$ such that $p_n < y < f^i(p_n)$ for all $1 \le i \le n-1$.  This confirms (b).

Here is a different proof of (b).  Let $m \ge 3$ be odd.  Let $\hat \ell$ be the {\it smallest odd} integer in $[1, m]$ such that $f^{\hat \ell}(\min P) < y$.  For $0 \le i \le \hat \ell -1$, let $I_i = [f^i(\min P) : f^i(y)]$.  If $\hat \ell = m$, then since $f^{\hat \ell-1}(\min P) = b$, we consider, for each $n \ge m+1$, the cycle $$I_0I_1I_2 \cdots I_{\hat \ell-2}([v, b])^{n-\hat \ell+1}I_0$$ of length $n$.  If $1 \le \hat \ell \le m-2$, then we consider, for each $n \ge m+1$, the cycle $$I_0I_1I_2 \cdots I_{\hat \ell-1} [v, f(y)] ([v, b])^{n-\hat \ell-1} I_0$$ of length $n$.  In either case, for each $n \ge m+1$, we have a period-$n$ point $q_n$ in $[\min P, y]$ such that $y < f^j(q_n)$ for all odd $1 \le j \le m$ and all $m+1 \le j \le n$.  This also confirms (b).  Note that the orbit of $q_n$ need not be the same as that of $p_n$ obtained above.

On the other hand, we can let $\tilde \ell$ be the {\it smallest odd} integer in $[1, m+2]$ such that $f^{\tilde \ell}(v) < y$.  For $0 \le i \le \tilde \ell -1$, let $J_i = [f^i(y) : f^i(v)]$.  If $\tilde \ell = m+2$, then $f^{\tilde \ell -1}(v) = f^{m+1}(v) = b$ and, for each $n \ge m+2$, we consider the cycle $$J_0J_1J_2 \cdots J_{\tilde \ell -2} ([v, b])^{n- \tilde \ell +1}J_0$$ of length $n$.  If $\tilde \ell \le m$, then, for each $n \ge m+1$, we consider the cycle $$J_0J_1J_2 \cdots J_{\tilde \ell -1} [v, f(y)] ([v, b])^{n- \tilde \ell -1} J_0$$ of length $n$.  In either case, we obtain, for each $n \ge m+2$, a period-$n$ point of $f$ whose orbit may be distinct from those of the points $p_n$ and $q_n$ obtained above.  This  confirms (b) too.  

For the proof of (c), let $m \ge 3$ be odd and let $z_0 = \min \{ v \le x \le z : f^2(x) = x \}$.  Then we have $f^2(x) < x < z_0 \le z < f(x)$ on $[v, z_0)$.  If $f^2(x) < z_0$ on $[\min P, v]$, then $f^2(x) < z_0$ on $[\min P, z_0)$.  Consequently, $(\min P \le) \, (f^2)^i(\min P) < z_0$ for each $i \ge 1$.  Since $m \ge 3$ is odd, this contradicts the fact that $(f^2)^{(m-1)/2}(\min P) = b > z_0$.  So, the point $d = \max \{ \min P \le x \le v : f^2(x) = z_0 \}$ exists and $f^2(x) < z_0 \le z < f(x)$ on $(d, v)$.  Let $I_0 = [d, v]$ and $I_1 = [v, z_0]$.  Then $f^2(I_0) \cap f^2(I_1) \supset I_0 \cup I_1$.  By considering the cycles $I_0 f(I_0) I_0$, $(I_0 f(I_0))^i I_1 f(I_1) I_0, i \ge 1$, we see that $f$ has periodic points of all {\it even} periods.  This establishes (c).

\section{The second digraph proof of (a), (b) and (c)}
Let $a$ and $b$ be points in $P$ such that $f(a) = \max P$ and $f(b) = \min P$.  If $b < a$, let $z$ be the {\it smallest} fixed point of $f$ in $(b, a)$.  Since $f^k(\min P) = \max P > z$ for some $k \ge 1$, we have $\max \{ f(x) : \min P \le x \le b \} > z$.  By considering the cycles $[b, z]$ $([\min P, b])^i$ $[b, z]$, $i \ge 0$, we obtain periodic points of all periods for $f$.  So, {\it for the rest of this section}, we assume that $a < b$.  

Let $v$ be the unique point in $P$ such that $f(v) = b$.  If $b < v$, then by considering the cycle $[b, v] ([a, b])^i [b, v]$, $i \ge 1$, we obtain that $f$ has periodic points of all periods $\ge 2$.  If $\min P < v < b$ and there is a fixed point $\tilde z$ of $f$ in $[\min P, v]$, then from the cycles $[\tilde z, v] ([v, b])^i [\tilde z, v]$, $i \ge 1$, we obtain that $f$ has periodic points of all periods $\ge 2$.  So, {\it for the rest of this section}, we assume that $\min P < v < b$ and $f$ has no fixed points in $[\min P, v]$.  Since $f^2$ is one-to-one on $P$, the set $f^2([\min P, v] \cap P)$ contains as many points as the set $[\min P, v] \cap P$.  Since $f^2(v) = \min P$, we have $f^2([\min P, v]) \supset [\min P, v]$.  So, there is a point $y$ in $[\min P, v]$ such that $f^2(y) = y$.  Since $f$ has no fixed points in $[\min P, v]$, $y$ is a period-2 point of $f$ such that $\min P < y < v < f(y)$.  This fact about $y$ will be used below.  This proves (a).  

We now give a proof of (b).  Let $m \ge 3$ be odd.  Since $f^{m+2}([y, v]) \supset [f^{m+2}(y) : f^{m+2}(v)] = [\min P, f(y)] \supset [y, v]$, there is a point $c$ in $[y, v]$ such that $f^{m+2}(c) = c$.  If $c$ has least period $m+2$, we are done.  Otherwise, let $n$ be the least period of $c$ under $f$.  Then since we have assumed that $f$ has no fixed points in $[\min P, v]$, $n \ge 3$ and is a proper divisor of $m+2$ which is odd and so $m+2 \ge 3n > 2n$.  Consequently, $m+2 - n > (m+2)/2$.  Since $[f^{n}(y) : f^{n}(c)] = [c, f(y)] \supset [v, f(y)]$ and $f([v, b]) \supset [\min P, b]$, by considering the cycle $$[y, c] [f(y) : f(c)] [f^2(y) : f^2(c)] \cdots [f^{n-1}(y) : f^{n-1}(c)] [v, f(y)] ([v, b])^{m+1-n} [y, c]$$ of length $m+2$, we obtain a point $d$ in $[y, c]$ such that $f^{m+2}(d) = d$, $f^{n}(d) \in [v, f(y)]$ and $f^i(d) \in [v, b]$ for all $n+1 \le i \le m+1$.  Thus, $y < d < v < f^j(d)$ for all $n \le j \le m+1$.  Since there are at least $(m+1)-n+1$ $(> (m+2)/2)$ {\it consecutive} iterates of $d$ under $f$ lying to the right of $v$ and $y < d < v$, the least period of $d$ under $f$ is $m+2$.  This confirms (b).

Now since $\min P = f(b) \le a < b \le \max P = f(a)$, there exist points $a < u < w < b$ such that $f(u) = b$ and $f(w) = a$.  Thus, if $m \ge 3$ is odd, then $P$ is also a period-$m$ orbit of $f^2$ and $\min P = f^2(u) < u < w < f^2(w) = \max P$.  By mimicking the proof in the first paragraph, we obtain that $f^2$ has a period-3 point.  So, $f$ has either a period-6 point or a period-3 point.  In either case, $f$ has a period-6 point.  Now that we have shown the existence of period-6 points for $f$, to complete the proof of (c), it suffices to prove the existence of period-$(2m)$ points for $f$.  For this purpose, we present three different approaches as follows:  

The first approach is independent of the existence of period-6 points of $f$ and is very similar to the one as in the proof of (b) above.  If $m = 3$, then $f$ has periodic points of all periods and we are done.  If $m \ge 5$ and odd, then since $f([f^m(y) : f^m(v)]) \supset [f^{m+1}(y) : f^{m+1}(v)] = [y, b] \supset [v, b]$, by considering the cycle $$[y, v] [f(y) : f(v)] [f^2(y) : f^2(v)] [f^3(y) : f^3(v)] \cdots [f^m(y) : f^m(v)] ([v, b])^{m-1} [y, v]$$ of length $2m$, we obtain a period-$(2m)$ point $p_{2m}$ of $f$ such that $\min P < f^2(p_{2m}) < y < p_{2m} < v < f^k(p_{2m}) < b$ for all $m+1 \le k \le 2m-1$.  

The second approach is based on the fact {\bf{\cite{al, bc, bh}}} that if $f$ has a period-6 point then $f$ has a period-6 point whose iterates are "jumping" alternately around a fixed point of $f$ (such orbits are called {\it simple} in the literature).  Indeed, if $\{ x_1, x_2, \cdots, x_6 \}$ with $x_1 < x_2 < \cdots < x_6$ is a period-6 orbit of $f$, let $s = \max \{ 1 \le i \le 6 : f(x_i) > x_i \}$.  Then $f(x_s) \ge x_{s+1}$ and $f(x_{s+1}) \le x_s$ and there is a fixed point $z$ of $f$ in $(x_s, x_{s+1})$.  If there is no $0 \le i \le 5$ such that $f^i(x_s)$ and $f^{i+1}(x_s)$ lie on the same side of $z$, then $x_s$ is such a period-6 point of $f$. Otherwise, let $r$ be the {\it smallest} integer in $[1, 5]$ such that $f^r(x_s)$ and $f^{r+1}(x_s)$ lie on the same side of $z$.  Then $1 \le r \le 4$.  We now consider appropriate cycles of length 6 as follows: \\
\indent If $r = 1$, the cycle \, $[x_s, z] [x_{s+1}, f(x_s)] ([x_s, z][z, x_{s+1}])^2 [x_s, z]$, \\  
\indent if $r = 2$, the cycle \, $[x_s, z] [z, f(x_s)] [f^2(x_s), x_s] ([z, x_{s+1}][x_s, z])^2$, \\
\indent if $r = 3$, the cycle \, $[x_s, z] [z, f(x_s)] [f^2(x_s), z] [x_{s+1}, f^3(x_s)] [x_s, z] [z, x_{s+1}] [x_s, z]$, \\
\indent if $r = 4$, the cycle \, $[x_s, z] [z, f(x_s)] [f^2(x_s), z] [z, f^3(x_s)] [f^4(x_s), x_s] [z, x_{s+1}] [x_s, z]$. \\  Or, we can choose points $x_s = x_0 < x_{-2} < x_{-4} < z < x_{-3} < x_{-1} < x_{s+1}$ such that $f(x_{-i}) = x_{-i+1}$, $1 \le i \le 4$, and consider the cycle $$[x_{-5+r} : z] [x_{-4+r} : z] [x_{-3+r} : z] \cdots [x_s : z] [f(x_s) : z] [f^2(x_s) : z] \cdots [f^{r-1}(x_s) : z] \, L \, [x_{-5+r} : z]$$ of length 6, where $L = [f^r(x_s), x_s]$ if $f^r(x_s) < z$, and $L = [x_{s+1}, f^r(x_s)]$ otherwise.  In either case, $f$ has a period-6 point whose iterates are jumping alternately around the fixed point $z$ of $f$.  Consequently, the left three points form a period-3 orbit for $f^2$.  By using the two adjcent compact intervals formed by these three period-3 points of $f^2$, we can find, for each $n \ge 1$, a period-$n$ orbit $Q_n$ of $f^2$ such that $Q_n \cup f(Q_n)$ is a period-$(2n)$ orbit of $f$.  So, if $f$ has a period-6 point then $f$ has periodic points of all even periods, including period $2m$.  

Here is the third approach: Suppose $f$ has a period-6 point.  Then $f^2$ has a period-3 orbit and so, for each {\it odd}  $n \ge 5$, $f^2$ has a \v Stefan cycle of least period $n$.  Indeed, without loss of generality, we may assume that $\{ a_1, a_2, a_3 \}$ is a period-3 orbit of $f^2$ with $f^2(a_3) = a_1 < a_2 = f^2(a_1) < a_3 = f^2(a_2)$.  Let $\hat z$ be a fixed point of $f^2$ in $(a_2, a_3)$.  Since $f^2(a_3) < a_2 < \hat z < a_3 = f^2(a_2)$, there are points $u_0, u_{-1}, u_{-2}, \cdots$ such that $f^2(u_{-j}) = u_{-j+1}$, $j \ge 1$ and $$a_2 = u_o < u_{-2} < \cdots < u_{-2i+2} < u_{-2i} < \cdots < \hat z < \cdots < u_{-2i+1} < u_{-2i+3} < \cdots < u_{-1} < a_3.$$  Let $J = [a_1, a_2]$.  For each $k \ge 2$, by considering the cycle (with respect to $f^2$) $$J [u_{-2k+2}, \hat z] [\hat z, u_{-2k+3}] [u_{-2k+4}, u_{-2k+2}] [u_{-2k+3}, u_{-2k+5}] \cdots [u_{-2}, u_{-4}] [u_{-3}, u_{-1}] [a_2, u_{-2}] [u_{-1}, a_3] J$$ of length $2k+1$, we obtain a point $w_k$ in $J = [a_1, a_2]$ such that $(f^2)^{2k+1}(w_k) = w_k$ and $$a_1 < w_k < a_2 < f^{4k-2}(w_k) < \cdots < f^6(w_k) < f^2(w_k) < \hat z < f^4(w_k) < \cdots < f^{4k}(w_k) < a_3.$$  So, $w_k$ is a periodic point of $f$ with least period $2k+1$ or $4k+2$.  If $w_k$ is a period-$(4k+2)$ point of $f$, then we are done.  Otherwise, we have $f^{2k+4}(w_k) = f^3(w_k)$ and $f^{2k+6}(w_k) = f^5(w_k)$.  Consequently, $f([f^2(w_k), f^4(w_k)]) \supset [f^3(w_k) : f^5(w_k)] = [f^{2k+4}(w_k): f^{2k+6}(w_k)] \supset [f^2(w_k), f^4(w_k)]$.  Therefore, $f$ has a fixed point in $[f^2(w_k), f^4(w_k)]$.  Without loss of generality, we may assume that $\hat z$ is a fixed point of $f$.  If $k$ is odd, then $f^3(w_k) = f^{2k+4}(w_k) < f^2(w_k) < \hat z < f^4(w_k)$ and we consider the cycle $[f^2(w_k), \hat z] [f^2(w_k), \hat z] [f^3(w_k), f^2(w_k)] [f^2(w_k), \hat z]$.  If $k$ is even, then $f^5(w_k) = f^{2k+6}(w_k) < f^2(w_k) < \hat z < f^4(w_k) < f^{2k+4}(w_k) = f^3(w_k)$ and we consider the cycle $[f^2(w_k), \hat z] [f^4(w_k), f^3(w_k)] [\hat z, f^4(w_k)] [f^2(w_k), \hat z]$.  In either case, we obtain a period-3 point of $f$ and so $f$ has periodic points of all periods, including a period-$(4k+2)$ point.  This shows that $f$ has a period-$(2j)$ point for each odd integer $j \ge 3$, including period $2m$.

For simplicity, in the following {\it five} sections, we let $P = \{x_i : 1 \le i \le m \}$, with $x_1 < x_2 < \cdots < x_m$, be a period-$m$ orbit of $f$ with $m \ge 3$ and let $s$ be an integer in $[1, m-1]$ such that $f(x_s) \ge x_{s+1}$ and $f(x_{s+1}) \le x_s$.  Let $z$ be a fixed point of $f$ in $(x_s, x_{s+1})$.  

\section{The third digraph proof of (a), (b) and (c)}
We first show by induction on $n \ge 2$ that if there exist a fixed point $z^*$ of $f$, a point $c$ and an integer $n \ge 2$ such that $f(c) < c < z^* < f^n(c)$ then there exist a fixed point $\tilde z$ of $f$ and a point $d$ such that $f(d) < d < \tilde z < f^2(d)$.  Indeed, if $n = 2$, we are done.  If $f^2(c) < f(c) < c < z^* < f^n(c)$, then $f(f(c)) < f(c) < z^* < f^{n-1}(f(c))$.  If $f(c) < f^2(c) < c < z^*$, let $\hat c$ be a point in $(c, z^*)$ such that $f(\hat c) = f^2(c)$, then $f(\hat c) < \hat c < z^* < f^n(c) = f^{n-1}(\hat c)$.  If $f(c) < c < f^2(c) < f^3(c)$, then there is a fixed point $\hat z$ of $f$ in $(c, f^2(c))$ and we are done.  If $f(c) < c < f^2(c) < z^*$ and $f^3(c) < f^2(c)$, let $\hat c$ be a point in $(f^2(c), z^*)$ such that $f(\hat c) = f^2(c)$, then we have $f(\hat c) < \hat c < z^* < f^n(c) = f^{n-1}(\hat c)$.  
In either case, we obtain by induction hypothesis a fixed point $\tilde z$ of $f$ and a point $d$ such that $f(d) < d < \tilde z < f^2(d)$.  Consequently, let $u$ be a point in $(f(d), d)$ such that $\tilde z < f(u) < f(d)$ and let $w$ be a point in $(d, \tilde z)$ such that $f(w) = u$.  Let $I_0 = [u, d]$ and $I_1 = [w, \tilde z]$.  Then $I_0 \cap I_1 = \emptyset$ and $f(I_0) \cap f(I_1) \supset I_0 \cup I_1$.  By considering the cycles $I_0(I_1)^nI_0$, $i \ge 1$, we obtain that $f$ has periodic points of all periods.  Now, let $a$ and $b$ be two points in $P$ such that $f(a) = x_m$ and $f(b) = x_1$.  If $b < x_s$ (if $a > x_{s+1}$, the argument is similar), then $x_1 = f(b) < b < x_s < z < x_m = f^k(b)$ for some $k \ge 2$.  The above implies that $f$ has periodic points of all periods and we are done.  So, {\it for the rest of this section}, we assume that $a \le x_s < z < x_{s+1} \le b$.

It is clear that one side of $z$ contains at least as many points of $P$ as the other side.  {\it We may assume that it is the right side}, (if it is the left side, the proof is similar).  If the right side contains as many points of $P$ as the left side, then, since $f([x_{s+1}, x_m])$ contains $f(b) = x_1$ and contains as many points of $P$ as $[x_1, x_s]$, we have $f([x_{s+1}, x_m]) \supset [x_1, x_s]$.  Similarly, $f([x_1, x_s]) \supset [x_{s+1}, x_m]$.  If the right side contains more points of $P$ than the left side, then we have $f([x_{s+1}, x_m]) \supset [x_1, x_{s+1}] \supset [x_1, z]$ and $f([x_1, z]) \supset f([a, z]) \supset [x_{s+1}, x_m]$.  In either case, $f$ has a period-2 point.  This proves (a).    

For the proof of (b) and (c), we assume that $m \ge 3$ is odd and the right side of $z$ contains more points of $P$ than the left side.  So, we have $1 \le s < m/2$ and $f([x_{s+1}, x_m]) \supset [x_1, x_{s+1}]$.  If $s = 1$, then $f(x_2) = f(x_{s+1}) = x_1$ and  $f(x_1) = f(a) = x_m$, and so $f([x_1, x_2]) \supset [x_1, x_m]$.  By considering the cycles $[x_2, x_m] ([x_1, x_2])^n [x_2, x_m]$, $n \ge 1$, we have periodic points of all periods $\ge 2$ for $f$ and we are done.  So, {\it for the rest of this section}, suppose $2 \le s < m/2$ (and so $m \ge 5)$.  We may also suppose $f$ has no period-$\ell$ points with $\ell$ odd and $1 < \ell < m$.    

Let $L = [x_{s+1}, x_m]$.  Then $f(L) \supset [x_1, x_{s+1}] \supset [a, x_{s+1}]$ and so $f^2(L) \supset [x_s, x_m] \supset L$.  Since $P$ is also a period-$m$ orbit of $f^2$, $f^2(P \cap [x_{s+1-k}, x_m]) \not\subset [x_{s+1-k}, x_m]$ for each $0 \le k \le s-1$.  Therefore, since $x_m \in f^2(L)$, we have $f^2([x_{s+1-k}, x_m]) \supset [x_{s-k}, x_m]$ for all $0 \le k \le s-1$.  Consequently, we have $f^{2s}(L) \supset [x_1, x_m]$.  By considering the cycle $L f(L) f^2(L) f^3(L) \cdots f^{2s-1}(L) [a, z] L$ of length $2s+1 \,(\le m)$, we obtain a period-$\ell$ point of $f$ with $\ell$ odd and $1 < \ell \le 2s+1 \le m$.  Since $f$ has no periodic points of odd periods $\ell$ with $1 < \ell < m$, this forces $2s + 1 = m$.  If $x_1 < a \le x_s$, or $f^2([x_{s+1-k}, x_m]) \supset [x_{s-k-1}, x_m]$ for some $0 \le k \le s-2$, then since $f^2([x_1, x_m]) \supset [x_1, x_m]$, we have $f^{2(s-1)}(L) \supset [x_1, x_m]$.  So, by considering the cycle $L f(L) f^2(L) f^3(L) \cdots f^{2s-3}(L) [a, z] L$ of length $2s-1 = m-2$, we obtain a period-$\ell$ point of $f$ with $\ell$ odd and $1 < \ell < m$.  This is a contradiction.  Therefore, $a = x_1$ and $f^2([x_{s+1-k}, x_m]) \not\supset [x_{s-k-1}, x_m]$ for all $0 \le k \le s-2$.  In particular, $x_{s-k-1} \notin f^2([x_{s+1-k}, x_m])$ for all $0 \le k \le s-2$.  Since $f^2([x_{s+1-k}, x_m]) \supset [x_{s-k}, x_m]$ for all $0 \le k \le s-1$, we obtain that $P \cap f^2([x_{s+1-k}, x_m]) \subset [x_{s-k}, x_m]$ for all $0 \le k \le s-1$.  

Since $x_{s+1} \le b \le x_m$ and since $P$ has more points on the right side of $z$ than on the left, we have $x_{s+1} \in [x_1, x_{s+1}] \subset f([x_{s+1}, x_m])$.  So, $f(x_{s+1}) \in P \cap f^2([x_{s+1}, x_m]) \subset [x_s, x_m]$.  In particular, $f(x_{s+1}) \ge x_s$.  But since $f(x_{s+1}) \le x_s$, this implies $f(x_{s+1}) = x_s$.  If $f^2(x_s) \in [x_s, x_m]$, then since $P \cap f^2([x_{s+1}, x_m]) \subset [x_s, x_m]$, we have $f^2(P \cap [x_s, x_m]) \subset \{ f^2(x_s) \} \cup \bigl(P \cap f^2([x_{s+1}, x_m])\bigr) \subset [x_s, x_m]$ which is a contradiction.  So, $f^2(x_s) < x_s$.  Since $P \cap f^2([x_s, x_m]) \subset [x_{s-1}, x_m]$, we have $f^2(x_s) = x_{s-1}$.  Inductively, we obtain that $f^2(x_{s+1-k}) = x_{s-k}$ for each $1 \le k \le s-1$.  In particular, $x_{s-k} = f^{2k}(x_s)$ for each $1 \le k \le s-1$ and $x_m = f(a) = f(x_1) = f(f^{2(s-1)}(x_s)) = f^{m-2}(x_s)$.  Since $f(x_{s+1}) = x_s$, this shows that $f^j(x_{s+1}) < z < f^i(x_{s+1}) < f^{m-1}(x_{s+1}) = x_m$ for all odd $1 \le j < m$ and all even $0 \le i < m-1$.  For each $n > m$, by considering the cycle $$[z, x_{s+1}] [z : f(x_{s+1})] [z : f^2(x_{s+1})] \cdots [z : f^{m-2}(x_{s+1})] [x_{m-1}, x_m] ([x_s, x_{s+1}])^{n-m} [z, x_{s+1}]$$ of length $n$, we obtain a period-$n$ point of $f$.  This establishes (b).  As for (c), it follows from Lemma 6 with $d = f^{m-3}(x_{s+1})$ (note that $f(x_m) = x_{s+1} > z$).  

\section{The fourth digraph proof of (a), (b) and (c)}
We now consider how the iterates $x_s, f(x_s), f^2(x_s), \cdots, f^{m-1}(x_s)$ "jump" around $z$.  

If for all integers $k$ such that $0 \le k \le m-1$ the points $f^k(x_s)$ and $f^{k+1}(x_s)$ lie on opposite sides of $z$, then, since $x_s < z < f(x_s)$, $m$ is even and, $f^i(x_s) < z < f^j(x_s)$ for all even $i$ and all odd $j$ in $[0, m-1]$.  This implies that $f([x_1, x_s] \cap P) = [x_{s+1}, x_m] \cap P$ and $f([x_{s+1}, x_m] \cap P) = [x_1, x_s] \cap P$.  In particular, $f([x_1, x_s]) \supset [x_{s+1}, x_m]$ and $f([x_{s+1}, x_m]) \supset [x_1, x_s]$.  By considering the cycle $[x_1, x_s][x_{s+1}, x_m][x_1, x_s]$, we obtain a period-2 point of $f$.  

We know the points $x_s$ and $f(x_s)$ always lie on the opposite sides of $z$.  So, assume that there is a {\it smallest} integer $r$ in $[1, m-1]$ such that the points $f^r(x_s)$ and $f^{r+1}(x_s)$ lie on the same side of $z$ (and so the iterates $x_s, f(x_s), f^2(x_s), \cdots, f^r(x_s)$ are jumping around $z$ alternately and this includes the case when $m$ is odd).  If both $f^r(x_s)$ and $f^{r+1}(x_s)$ lie on the right side of $z$, then $r$ is odd and we consider the cycle $$[x_s, z][z, f(x_s)][f^2(x_s), z] \cdots [f^{r-1}(x_s), z][x_{s+1}, f^r(x_s)] [x_s, z].$$If both $f^r(x_s)$ and $f^{r+1}(x_s)$ lie on the left side of $z$, then $r$ is even and we consider the cycle $$[x_s, z][z, f(x_s)] [f^2(x_s), z] \cdots [z, f^{r-1}(x_s)][f^r(x_s), x_s][z, x_{s+1}][x_s, z].$$  In either case, we obtain a periodic point of $f$ whose iterates are alternating around $z$.  It follows from the previous paragraph that $f$ has a period-2 point of $f$.  This, combined with the previous paragraph, implies (a).

Now suppose there is a {\it smallest} integer $1 \le r \le m-1$ such that both $f^r(x_s)$ and $f^{r+1}(x_s)$ lie on the same side of $z$.  If $r = 1$, we consider the cycles $[x_s, x_{s+1}] [x_{s+1}, f(x_s)] ([x_s, x_{s+1}])^i$, $i \ge 1$.  If $r = 2$, we consider the cycles $[x_s, x_{s+1}]$ $[z, f(x_s)]$ $[f^2(x_s), x_s]$ $([x_s, x_{s+1}])^i$, $i \ge 1$.  In either case, we have periodic points of all periods $\ge 3$, including period-$(m+2)$.  If $r \ge 3$, let $L_i = [f^i(x_s) : f^{i+2}(x_s)]$ for $0 \le i \le r-3$.  For each integer  $n \ge m+1$, we consider the cycle $$[x_s, x_{s+1}] [z, f(x_s)] L_0 L_1 L_2 L_3 \cdots L_{r-3} [f^{r-2}(x_s) : f^r(x_s)]([x_s, x_{s+1}])^{n-r}$$ of length $n$ to obtain a period-$n$ point $p_{n}$ of $f$ in $[x_s, x_{s+1}]$ such that $f^i(p_n) \not\in [x_s, x_{s+1}]$ for all $2 \le i \le r$ and $f^j(p_n) \in [x_s, x_{s+1}]$ for all $r+1 \le j \le n$.  In particular, this shows that if $f$ has a period-$m$ point with $m \ge 3$ and odd, then $f$ has a period-$(m+2)$ point.  That is, (b) holds.  

On the other hand, let $1 \le k \le r \, (\le m-1)$ be the smallest integer such that $[z : f^k(x_s)] \supset [z : f^r(x_s)]$.  If $k = 1$, then by considering the cycles $[x_s, x_{s+1}] [x_{s+1}, f(x_s)] ([x_s, x_{s+1}])^i$, $i \ge 1$, we have periodic points of all periods $\ge 2$ and we are done.  So, suppose $k \ge 2$.  Then $[z : f^{k-2}(x_s)] \subsetneq [z : f^r(x_s)]$.  By considering the cycles $$[f^{k-2}(x_s): f^r(x_s)] [z : f^{k-1}(x_s)] \,\, ([z : f^{k-2}(x_s)] [z : f^{k-1}(x_s)])^i \,\, [f^{k-2}(x_s): f^r(x_s)], \, i \ge 0,$$ we obtain periodic points of all even periods $\ge 2$, including period-6 and period-$(2m)$ for $f$.  That is, (c) holds.  

In the following, we present different proofs of (b) and (c).

We may assume that, {\it for the rest of this section}, both $f^r(x_s)$ and $f^{r+1}(x_s)$ lie on the right side of $z$ (if they both lie on the left side of $z$, the proof is similar).  Then $r$ is odd and $1 \le r \le m-2$.  Since $f(x_{s+1}) \le x_s < x_{s+1}$, we have $f^r(x_s) > x_{s+1}$.  Let $k$ be the {\it smallest} integer in $[1, r]$ such that $f^k(x_s) \ge f^r(x_s)$.  Since the iterates $x_s, f(x_s), f^2(x_s), \cdots$, $f^k(x_s)$ are "jumping" alternately around $z$, we see that $k$ is odd.  If $k = 1$, then for all $i \ge 2$, by considering the cycles $[x_{s+1}, f^r(x_s)]([x_s, x_{s+1}])^{i-1}[x_{s+1}, f^r(x_s)]$ of length $i$, we see that $f$ has periodic points of all periods $\ge 2$.  So, suppose $k \ge 3$.  Thus, $f^{k-1}(x_s) < z < f^{k-2}(x_s) < f^r(x_s) \le f^k(x_s)$.  Let $U = [f^{k-1}(x_s), z]$, $V = [z, f^{k-2}(x_s)]$ and $W = [f^{k-2}(x_s), f^r(x_s)]$.  Then $f(U) \supset V \cup W$, $f(V) \supset U$ and $f(W) \supset U$.  By considering the cycles $WUW$ and $W(UV)^iUW, i \ge 1$, we establish (c) and also (a).  

On the other hand, let $x_0, x_{-1}, x_{-2}, \cdots$ be points in $[x_s, x_{s+1}]$ such that $$x_s = x_0 < x_{-2} < x_{-4} < \cdots < z < \cdots < x_{-3} < x_{-1} < x_{s+1}$$ and $f(x_{-i}) = x_{-i+1}$ for all $i = 1, 2, \cdots$.  For each {\it odd} integer $n \ge m$, let $n_0 = n-1-r$.  Note that we suppose $r$ is odd and so $n_0$ is odd.  By considering the cycle $$[f^{r-2}(x_s) : f^{r}(x_s)] [z, x_{-n_0}] [x_{-n_0+1}, z] [x_{-n_0}, x_{-n_0+2}] [x_{-n_0+3}, x_{-n_0+1}] \cdots [x_{-3}, x_{-1}] [x_s, x_{-2}] \quad\,\,\,\,$$ $$\qquad\qquad\qquad\qquad [x_{-1}, f(x_s)] [f^2(x_s), x_s] [f(x_s) : f^3(x_s)] [f^2(x_s) : f^4(x_s)] \cdots [f^{r-2}(x_s) : f^{r}(x_s)]$$ of length $n$, we obtain a period-$n$ point $p_n$ in $[f^{r-2}(x_s) : f^r(x_s)]$ whose iterates $f^i(p_n)$, $1 \le i \le n$, are jumping alternately around $z$.  This confirms (b).  We can also consider the cycle $$[x_{s+1}, f^r(x_s)] ([z, x_{s+1}][x_s, z])^{(n-r)/2} [z, f(x_s)] [z : f^2(x_s)] \cdots [z : f^{r-1}(x_s)] [x_{s+1}, f^r(x_s)]$$ of length $n$ to obtain a period-$n$ point $q_n$ in $[x_{s+1}, f^r(x_s)]$ such that the iterates $f^i(q_n)$, $1 \le i \le n$, are jumping alternately around $z$.  This also proves (b).  Note that the orbit of $q_n$ need not be the same as that of the above point $p_n$.  

With a little more work than the above proof of (b), we can even show the existence of all \v Stefan cycles of odd periods $n \ge m$ {\it without assuming the non-existence of odd periods $\ell$ with $1 < \ell < m$}.  Indeed, since we suppose $r$ is odd, $1 \le r \le m-2$.  Let $u$ be a point in $[z, x_{s+1}]$ such that $f(u) = x_s$.  By considering the cycle $$[z, u] [x_s, z] [z, f(x_s)] [f^2(x_s), z] [z, f^3(x_s)] \cdots [f^{r-1}(x_s), z] [x_{s+1}, f^r(x_s)] [z, u]$$ of length $r+2 \le m$, we obtain a period-$(r+2)$ point $w$ of $f$ in $[z, u]$ such that the iterates $f^j(w)$, $0 \le j \le r+1$, are jumping alternately around $z$ and $z < w < f^{r+1}(w)$.  For $0 \le i \le r+1$, let $J_{i} = [z : f^{i}(w)]$.  If the orbit of $w$ with respect to $f$ is not a \v Stefan cycle, then, for some $0 \le i_2 < i_1 \le r+1$ with $i_1 - i_2$ even, we have $[f^{i_1}(w) : z] \subset [f^{i_2}(w) : z]$.  If $i_2 = 0$, then $i_1 < r+1$ (since $z < w < f^{r+1}(w)$) and we consider the cycle $J_{i_1}J_{i_1+1}J_{i_1+2} \cdots J_r[w, f^{r+1}(w)]J_{i_1}$ of length $r+2 - i_1$.  Otherwise, we consider the cycle $J_0J_1 \cdots J_{i_2-1}J_{i_1}J_{i_1+1} \cdots$ $J_{r}$ $[w, f^{r+1}(w)]J_0$ if $i_1 < r+1$ or the cycle $J_0J_1J_2 \cdots J_{i_2-1} [w, f^{r+1}(w)]J_0$ if $i_1 = r+1$ of odd length $r+2 - i_1 + i_2$.  In either case, we obtain a periodic point $\hat w$ of $f$ with {\it smaller} odd period $\hat r +2$, where $\hat r = r -i_1+i_2$ and $1 < \hat r +2 < r+2$ such that the iterates $f^j(\hat w)$, $0 \le j \le \hat r +1$, are jumping alternately around $z$ and $z < \hat w < f^{\hat r +1}(\hat w)$.  Proceeding in this manner finitely many times, we eventually obtain a periodic point $q$ of odd period $1 < \ell < m$ such that $$f^{\ell-2}(q) < f^{\ell-4}(q) < \cdots < f^3(q) < f(q) < z < q < f^2(q) < f^4(q) <\cdots < f^{\ell-3}(q) < f^{\ell-1}(q).$$That is, The orbit of $q$ with respect to $f$ is a \v Stefan cycle.  It is now easy to see that $f$ has, for each odd $n \ge \ell$ (and so for each $n \ge m$), a \v Stefan cycle of odd period $n$.  

\section{The fifth digraph proof of (a), (b) and (c)}
We now reconsider how the iterates $x_s$, $f(x_s)$, $f^2(x_s)$, $\cdots$, $f^{m-1}(x_s)$ "jump" around the fixed point $z$ under the assumption that, when $m \ge 3$ is odd, $f$ has no periodic points of smaller odd periods other than fixed points.

If for all integers $k$ such that $0 \le k \le m-1$ the points $f^k(x_s)$ and $f^{k+1}(x_s)$ lie on opposite sides of $z$, then, since $x_s < z < f(x_s)$, $m$ is even and, $f^i(x_s) < z < f^j(x_s)$ for all even $i$ and all odd $j$ in $[0, m-1]$.  This implies that $f([x_1, x_s]) \supset [x_{s+1}, x_m]$ and $f([x_{s+1}, x_m]) \supset [x_1, x_s]$.  Consequently, $f$ has a period-2 point.  

Now assume that there is a {\it smallest} integer $r$ in $[1, m-1]$ such that the points $f^r(x_s)$ and $f^{r+1}(x_s)$ lie on the same side of $z$ (and so the iterates $x_s$, $f(x_s)$, $f^2(x_s)$, $\cdots, f^r(x_s)$ are jumping around $z$ alternately and this includes the case when $m$ is odd).  The following proof of (a) is not needed here.  We include it for the sake of interest in itself.  For any finite set $A$ of real numbers, let $H(A)$ denote the interval $[\min A, \max A]$.  If $r$ is even, let $I_0 = H(\{ x_s, f^2(x_s), f^4(x_s), \cdots, f^r(x_s)\})$ and $I_1 = H(\{z, f(x_s), f^3(x_s), \cdots, f^{r-1}(x_s)\})$.  If $r$ is odd, let $I_0 = H(\{ x_s, f^2(x_s), f^4(x_s), \cdots, f^{r-1}(x_s)$, $z\})$ and $I_1 = H(\{x_{s+1}, f(x_s), f^3(x_s), \cdots$, $f^r(x_s)\}) \, (x_{s+1}$ is needed only when $r = 1$).  In either case, $I_0 \cap I_1 = \emptyset$ and, $f(I_0) \supset I_1$ and $f(I_1) \supset I_0$.  So, $f$ has a period-2 point.  This, together with the above, confirms (a).

For the proofs of (b) and (c), let $m \ge 3$ be odd.  If $m = 3$, then $f$ has periodic points of all periods and we are done.  So, suppose $m > 3$.  Without loss of generality, we may also suppose $f$ has no non-fixed periodic points of smaller odd periods.  Since $m$ is odd, there is a {\it smallest} integer $1 \le r \le m-1$ such that the points $f^r(x_s)$ and $f^{r+1}(x_s)$ lie on the same side of $z$.  If $1 \le r \le m-3$, let $J_i = [z : f^i(x_s)]$ for all $0 \le i \le r-1$, and $J_r = [f^r(x_s), x_s]$ if $f^r(x_s)< x_s$ and $J_r = [x_{s+1}, f^r(x_s)]$ if $x_{s+1} < f^r(x_s)$.  By considering the cycle $J_0J_1J_2 \cdots J_r([x_s, x_{s+1}])^{m-3-r}J_0$ of length $m-2$, we obtain a non-fixed periodic point of $f$ with odd period $\le m-2$.  This is a contradiction.  So, for all $1 \le k \le m-3$, the points $f^k(x_s)$ and $f^{k+1}(x_s)$ lie on opposite sides of $z$.  Since $f(x_s) > z$, this implies that $f^i(x_s) < z < f^j(x_s)$ for all even $0 \le i \le m-3$ and all odd $1 \le j \le m-2$.  In particular, $z < f^{m-2}(x_s)$.  We now have three cases to consider depending on the locations of the point $f^{m-1}(x_s)$:

Case 1.  $f^{m-1}(x_s) < z < f^{m-2}(x_s)$.  In this case, we actually have $f^i(x_s) < z < f^j(x_s)$ for all even $0 \le i \le m-1$ and all odd $1 \le j \le m-2$.  For $0 \le i \le m-1$, let $J_i = [z : f^i(x_s)]$.  If, for some $1 \le \ell < k \le m-1$ with $k- \ell$ even, we have $[z : f^k(x_s)] \subset [z : f^\ell(x_s)]$, then by considering the cycle $J_0J_1 \cdots J_{\ell-1}J_kJ_{k+1} \cdots$ $J_{m-2}$ $[f^{m-1}(x_s), x_s]J_0$ if $k < m-1$ or the cycle $J_0J_1J_2 \cdots J_{\ell-1} [f^{m-1}(x_s), x_s]J_0$ if $k = m-1$, of odd length $m - k + \ell$, we obtain a periodic point of $f$ with odd period $m-k + \ell \le m-2$ other than fixed points.  This contradicts the assumption.  So, we must have 
$$
f^{m-1}(x_s) < \cdots < f^4(x_s) < f^2(x_s) < x_s < z < f(x_s) < f^3(x_s) < \cdots < f^{m-2}(x_s).
$$
\noindent
That is, $P$ is a \v Stefan cycle.  Consequently, for each $n > m$, by considering the cycle $$J_0J_1J_2 \cdots J_{m-2} [f^{m-1}(x_s), f^{m-3}(x_s)] ([x_s, f(x_s)])^{n-m} J_0$$ of length $n$, we obtain a period-$n$ point for $f$.  Furthermore, by Lemma 6 with $d = f^{m-3}(x_s)$, $f$ has periodic points of all even periods $\ge 2$.          

Case 2.  $x_{s+1} = f^{m-1}(x_s) < f^{m-2}(x_s)$.  In this case, we have $f^i(x_s) < z < x_{s+1} = f^{m-1}(x_s) < f^j(x_s)$ for all even $0 \le i \le m-3$ and all odd $1 \le j \le m-2$.  Since $x_{s+1} = f^{m-1}(x_s)$, we have $f(x_{s+1}) = x_s$.  Thus, $x_s = f(x_{s+1})$.  By plugging this in the above inequalities, we obtain that $f^j(x_{s+1}) < z < f^i(x_{s+1})$ for all odd $1 \le j \le m-2$ and all even $0 \le i \le m-1$.  This is a symmetric copy of Case 1.  Therefore, $P$ is a \v Stefan cycle and $f$ has periodic points of all periods $> m$ and of all even periods $\ge 2$.    

Case 3.  $x_{s+1} < f^{m-1}(x_s)$.  In this case, since $f(x_{s+1}) \le x_s < x_{s+1}$ and $f(f^{m-2}(x_s)) = f^{m-1}(x_s) > x_{s+1}$, we have $f^{m-2}(x_s) \ne x_{s+1}$.  Hence $x_{s+1} = f^k(x_s)$ for some $1 \le k \le m-3$.  So, $f^{m-k-2}(x_{s+1}) = f^{m-2}(x_s)$.  By considering the cycle $$[z, x_{s+1}] [z : f(x_{s+1})] [z : f^2(x_{s+1})] \cdots [z : f^{m-k-3}(x_{s+1})] [x_{s+1}, f^{m-2}(x_s)] ([x_s, x_{s+1}])^{k-1} [z, x_{s+1}]$$ of length $m-2$, we obtain a periodic point of $f$ with smaller odd period other than fixed points.  This contradicts the assumption.  So, this case cannot occur.  

\noindent
This shows that $P$ is a \v Stefan cycle and $f$ has periodic points of all periods $> m$ and of all even periods $\ge 2$.  
Therefore, (a), (b) and (c) are established.

\section{The sixth digraph proof of (a), (b) and (c)}
We now consider how points of $P$ which lie on either side of $z$ are mapped.  

If for all integers $i$ such that $1 \le i \le m-1$ and $i \ne s$ the points $f(x_i)$ and $f(x_{i+1})$ lie on the same side of $z$, then, since $f(x_{s+1}) < z < f(x_s)$, $f(P \cap [x_1, x_s]) \subset P \cap [x_{s+1}, x_m]$ and $f(P \cap [x_{s+1}, x_m]) \subset P \cap [x_1, x_s]$.  Since $f$ is one-to-one on $P$, we obtain that $f(P \cap [x_1, x_s]) = P \cap [x_{s+1}, x_m]$ and $f(P \cap [x_{s+1}, x_m]) = P \cap [x_1, x_s]$ (and so $m$ is even).  In particular, $f([x_1, x_s]) \supset [x_{s+1}, x_m]$ and $f([x_{s+1}, x_m]) \supset [x_1, x_s]$.  Consequently, $f$ has a period-2 point.    

Now {\it for the rest of this section}, assume that for some integer $t$ such that $1 \le t \le m-1$ and $t \ne s$ the points $f(x_t)$ and $f(x_{t+1})$ lie on opposite sides of $z$ (this includes the case when $m$ is odd).  Suppose $x_t < x_s$.  If $x_{s+1} \le x_t$, the proof is similar.  Since $x_t < x_s$ and $f(x_s) \ge x_{s+1}$, we may also assume that $t$ is the {\it largest} integer in $[1, s-1]$ such that $f(x_t) \le x_s$.  So, $f([x_t, x_{t+1}]) \supset [x_s, x_{s+1}]$ and $f(x) > z$ for all $x$ in $P \cap [x_{t+1}, z]$.  

If $f(x_i) \ge x_{t+1}$ for all $s+1 \le i \le m$, then $f^n(x_s) \ge x_{t+1}$ for all $n \ge 1$, contradicting the fact that $f^j(x_s) = x_t$ for some $1 \le j \le m-1$.  So, there is a {\it smallest} integer $\ell$ with $s+1 \le \ell \le m$ such that $f(x_\ell) \le x_t$.  If $x_{t+1} \le f(x_i) \le x_{\ell-1}$ for all $t+1 \le i \le \ell -1$, then $f^n(x_s) \ge x_{t+1}$ for all $n \ge 1$, again contradicting the fact that $f^j(x_s) = x_t$ for some $1 \le j \le m-1$.  So, there is an integer $k$ with $t+1 \le k \le \ell -1$ such that $f(x_k) \ge x_\ell$.  If $s+1 < k \le \ell -1$, by considering the cycles $[x_{s+1}, x_k]([x_k, x_\ell])^n[x_{s+1}, x_k]$, $n \ge 1$, we obtain periodic points of all periods $\ge 2$ for $f$ and we are done.  

{\it For the rest of this section}, suppose $t+1 \le k \le s$.  So, we have $x_t < x_k < z < z_\ell$.  Let $d$ be a point in $[x_k, z)$ such that $f(d) \in (z, x_\ell], f^2(d) \in [x_t, x_k]$ and $f^3(d) = z$.  By Lemma 6, we confirm the existence of periodic points of all even periods, including a period-2 point $y$ of $f$ in $[x_t, x_k]$ such that $x_t < y < x_k < z < f(y) < x_\ell$.  This, combined with the above, establishes (a).  Furthermore, this also shows that if $m \ge 3$ is odd then (c) holds.  

As for the proof of (b), let $m \ge 3$ be odd and let $x_t, y, x_k, z, x_\ell$ be defined as above.  Let $r$ be the {\it smallest odd} integer in $[2, m]$ such that $f^r(x_k) \le x_k$.  For each {\it odd} integer  $n > m$, by considering the cycle $$[x_t, x_k] [z, f(x_k)] [z : f^2(x_k)] \cdots [z : f^{r-1}(x_k)] ([x_k, z] [z, x_\ell])^{(n-r)/2} [x_t, x_k]$$ or the cycle $$[y, x_k] [f(y) : f(x_k)] [f^2(y) : f^2(x_k)] \cdots [f^{r-1}(y) : f^{r-1}(x_k)] ([x_k, z] [z, x_\ell])^{(n-r)/2} [y, x_k]$$ of length $n$, we obtain (from each cycle) a period-$n$ point of $f$ in $[x_t, x_k]$ whose odd ($< n+1$) iterates visit the interval $[\min I, x_k] \,\, (\supset [x_t, x_k] \supset [y, x_k])$ exactly once.  This proves (b).  Note that, so far we have shown that if there exists a point $p$ in $P$ such that $p$ and $f(p)$ lie on the same side of $z$ then $f$ has a period-$n$ point for each $n \prec m$ in the Sharkovsky ordering.    

We now consider iterates of the (different) point $x_\ell$ : Let $\hat r$ be the {\it smallest even} integer in $[2, m+1]$ such that $f^{\hat r}(x_\ell) \le x_k$.  For each {\it odd} integer  $n > m$, by considering the cycle $$[f(x_\ell), y] [f^2(x_\ell):f(y)] \cdots [f^{\hat r-1}(x_\ell) : f^{\hat r -2}(y)] ([x_k, z][z, x_\ell])^{(n- \hat r +1)/2} [f(x_\ell), y]$$ of length $n$, we obtain a period-$n$ point of $f$ in $[f(x_\ell), y] \,\,( \subset [\min I, x_k])$ whose {\it odd} ($< n+1$) iterates visit the interval $[\min I, x_k]$ exactly once.  This proves (b).  Note that the period-$n$ points obtained from the above three cycles may generate three distinct orbits, for example, if $m = 7, f^7(x_k) = x_k$ and $$x_t = f^6(x_k) < y < f^2(x_k) < x_k < z < f(y) < f^5(x_k) = x_\ell < f^3(x_k) < f(x_k) < f^4(x_k).$$   

If $f$ has a period-3 point, then it is clear that $f$ has \v Stefan cycles of all odd periods $\ge 3$.  So, in the following, we suppose $f$ has no period-3 points.  When $m \ge 5$ is odd, our method can even be used to prove the existence of \v Stefan cycles of all odd periods $n \ge m$ {\it without assuming the non-existence of periodic points of odd periods $< m$ other than fixed points}.  Indeed, as before, let $t$ be the largest integer in $[1, s-1]$ such that $f(x_t) \le x_s$.  Then $f(x_i) > z$ for all $t+1 \le i \le s$.  If $f(x_i) < z$ for all $s+1 \le i \le m$, then there is a {\it smallest} integer $\ell$ in $[s+1, m]$ such that $f(x_\ell) \le x_t$.  On the other hand,  if there is also a {\it smallest} integer $r$ in $(s+1, m]$ such that $f(x_r) > z$, then since $x_t \notin$ $[x_{t+1}, x_{r-1}]$, we have $f(P \cap [x_{t+1}, x_{r-1}]) \not\subset [x_{t+1}, x_{r-1}]$.  So, for some $\hat \ell$ in $[t+1, r-1]$, $f(x_{\hat\ell}) \not\in [x_{t+1}, x_{r-1}]$.  Thus, we have either $t+1 \le \hat\ell \le s$ and $f(x_{\hat\ell}) \ge x_r$ or, $s+1 \le \hat\ell < r$ and $f(x_{\hat\ell}) \le x_t$.  Since the proofs for these two cases and the above one are similar, we may assume that there is a {\it smallest} integer $\ell$ in $[s+1, m]$ such that $f(x_\ell) \le x_t$ and,   $f(x_i) > z$ for all $t+1 \le i \le s$ and $f(x_j) < z$ for all $s+1 \le j \le \ell - 1$.  Consequently, $f(x_k) \ge x_\ell$ for some $t+1 \le k \le s$.  

Since $f(x_k) \ge x_\ell$, the point $u_0 = \max \{ x_k \le x < z : f(x) = x_\ell \}$ exists and $f(x) < x_\ell$ for all $x$ in $(u_0, z)$.  If $P \cap (u_0, z) = \emptyset$, then $x_s \le u_0$.  In this case, let $j = 0$.  If $P \cap (u_0, z) \ne \emptyset$, let $p$ be a point in $P \cap (u_0, z)$.  Then $z < f(p) \le x_{\ell -1}$. It is clear that $f(p_{-1}) \le u_0$ for some $p_{-1}$ in $P \cap (z, x_{\ell-1}]$.  So, the point $u_{-1} = \min \{ z < x \le p_{-1} : f(x) = u_0 \}$ exists and $u_0 < f(x)$ for all $x$ in $(z, u_{-1})$.  Since $p \in [u_0, z]$, it is clear that $f(p_{-2}) \ge u_{-1}$ for some $p_{-2}$ in $P \cap (u_0, z)$.  So, the point $u_{-2} = \max \{ p_{-2} \le x < z : f(x) = u_{-1} \}$ exists and $f(x) < u_{-1}$ for all $x$ in $(u_{-2}, z)$.  If $P \cap (u_{-2}, z) = \emptyset$, then $x_s \le u_{-2}$.  Otherwise let $\tilde p$ be a point in $P \cap (u_{-2}, z)$.  Then $z < f(\tilde p) < u_{-1}$.  We proceed in this manner.  Since $(x_k, x_\ell)$ contains at most $m-3$ points of $P$ (exclude at least these three points $x_t, x_k, x_\ell$), there is an {\it even} integer $0 \le j \le m-3$ and points $p_{-1}, p_{-2}, \cdots, p_{-j}$ in $P$ and points $u_0$, $u_{-1}$, $u_{-2}, \cdots, u_{-j}$ in $[x_k, x_\ell]$ such that $P \cap (u_{-j}, z) = \emptyset$, $f(u_{-i}) = u_{-i+1}$ for all $1 \le i \le j$, and 
$$ 
x_k \le u_0 < p_{-2} \le u_{-2} < \cdots <  p_{-j+2} \le u_{-j+2} < p_{-j} \le u_{-j} < z \qquad\qquad\qquad\qquad\qquad \qquad
$$
$$
\qquad\qquad\qquad\qquad\qquad\qquad\qquad < u_{-j+1} \le p_{-j+1} < u_{-j+3} \le p_{-j+3} < \cdots < u_{-1} \le p_{-1} < x_\ell.
$$
It is clear that there exist more points $u_{-\hat i}$, $\hat i > j$, such that $f(u_{-\hat i}) = u_{-\hat i+1}$, $\hat i > j$ and $$u_{-j} < u_{-j-2} < u_{-j-4} < \cdots < z < \cdots < u_{-j-3} < u_{-j-1} < u_{-j+1}.$$  Since $P \cap (u_{-j}, z) = \emptyset$, we have $x_s \le u_{-j}$.  Since $f(x_t) \le x_s$,we have $f([x_t, x_k]) \supset [x_s, x_\ell] \supset [x_{-j}, z]$.  For each {\it even} integer $n \ge j$, by considering the cycle $$[x_t, x_k] [u_{-n}, z] [z, u_{-n+1}] [u_{-n+2}, u_{-n}] [u_{-n+1}, u_{-n+3}] \cdots [u_0, u_{-2}] [u_{-1}, x_\ell] [x_t, x_k]$$ of length $n+3 \,\, (\ge j+3)$, we obtain a \v Stefan cycle of least period $n+3$.  Since $j$ is even and $0 \le j \le m-3$, this shows that $f$ has a \v Stefan cycle of least period $n$  for each odd $n \ge m$.

If $f$ has no periodic points of odd periods in $[3, m-2]$, then $j$ must be equal to $m-3$ and there is exactly one point of $P$ in each of the $m-3$ half-open intervals $(u_0, u_{-2}]$, $(u_{-2}, u_{-4}]$, $\cdots, (u_{-j+2}, u_{-j}]$, $[u_{-j+1}, u_{-j+3})$, $\cdots, [u_{-5}, u_{-3})$, $[u_{-3}, u_{-1})$, $[u_{-1}, x_\ell)$.  These $m-3$ points, plus $x_k, x_\ell$, and $x_t$, constitute the orbit $P$.  Therefore, $P$ is itself a \v Stefan cycle.  

\section{The seventh digraph proof of (a), (b) and (c)}
We now reconsider (cf. {\bf{\cite{du2}}}) how the points of $P$ which lie on either side of the fixed point $z$ are mapped {\it under the assumption that, when $m \ge 5$ is odd, $f$ has no periodic points of smaller odd periods other than fixed points.}

If for all integers $i$ such that $1 \le i \le m-1$ and $i \ne s$ the points $f(x_i)$ and $f(x_{i+1})$ lie on the same side of $z$, then it is easy to see that $f([x_1, x_s]) \supset [x_{s+1}, x_m]$ and $f([x_{s+1}, x_m]) \supset [x_1, x_s]$.  So, $f$ has a period-2 point.  On the other hand, assume that there is an integer $1 \le t \le m-1$ such that $t \ne s$ and the points $f(x_t)$ and $f(x_{t+1})$ lie on opposite sides of $z$ (this includes the case when $m \ge 3$ is odd).  Without loss of generality, we may assume that $x_t < x_s$.  We may also assume that $t$ is the {\it largest} integer in $[1, s)$ such that $f(x_t) \le x_s$.  So, $f(x_i) \ge x_{s+1}$ for all $t+1 \le i \le s$.  

Now let $q$ be the {\it smallest} positive integer such that $f^q(x_s) \le x_t$.  Then $2 \le q \le m-1$.  If $q \le m-2$, then we consider the cycles $$[x_s, x_{s+1}] [z : f(x_s)] [z : f^2(x_s)] \cdots [z : f^{q-1}(x_s)] [x_t, x_{t+1}] ([x_s, x_{s+1}])^{i+1}, i \ge 0$$to obtain a periodic point of least period $q +i+1$ for each $i \ge 0$.  If $q = m-1$, we write $x_{s+1} = f^j(x_s)$ for some $1 \le j \le m-2$ and consider the cycles $$[z, f^j(x_s)] [z : f^{j+1}(x_s)] [z : f^{j+2}(x_s)] \cdots [z : f^{m-2}(x_s)] [x_t, x_{t+1}] ([x_s, x_{s+1}])^i [z, f^j(x_s)], i \ge 0$$to obtain a periodic point of least period $m-j+i$ for each $i \ge 0$.  In either case, we obtain periodic points of all periods $\ge m-1$.  Consequently, if $m \ge 3$ is odd, then this confirms (b), and if $m \ge 4$ is even, then $f$ has periodic points of all odd periods $\ge m-1$.  

For the proof of (c), we may assume that $m \ge 5$ is odd and $f$ has no periodic points of smaller odd periods other than fixed points.  Under this assumption, we want to show that $$f^{m-1}(x_s) = x_t < x_{t+1} \le f^{m-3}(x_s) < z < f^{m-2}(x_s).$$  Then by Lemma 6, we easily obtain (c) and (a).  Since we have assumed that $f$ has no periodic points of smaller odd periods other than fixed points, the above shows that $q$ cannot be $\le m-3$.  So, either $q = m-2$ or $q = m-1$, and $x_{t+1} \le f^i(x_s)$ for all $0 \le i \le m-3$.  We also have $f(x_i) \ge x_{s+1}$ for all $t+1 \le i \le s$.  Based on these facts, we present three different arguments below.  

Firstly, suppose $q = m-2$ and so $f^{m-2}(x_s) \le x_t$.  If $f^{m-1}(x_s) \ge x_{t+1}$, then $f^{m-1}(x_s) \ge x_{s+1}$.  But then $f^{m-2}(x_s) \ne x_t$ because $f(x_t) \le x_s$.  So, we have $x_t = f^k(x_s)$ for some $1 \le k \le m-3$.  By considering the cycle $$[x_s, x_{s+1}] [z, f(x_s)] [z : f^2(x_s)] \cdots [z : f^{k-1}(x_s)] [x_t, x_s] ([x_s, x_{s+1}])^{m-2-k}$$ of length $m-2$, we obtain a periodic point of samller odd period other than fixed point which is a contradiction.  So, if $f^{m-2}(x_s) \le x_t$, then we also have $f^{m-1}(x_s) \le x_t$.  But then $x_{s+1} = f^k(x_s)$ for some $1 \le k \le m-3$.  By considering the cycle $$[x_s, f^k(x_s)] [z : f^{k+1}(x_s)] [z : f^{k+2}(x_s)] \cdots [z : f^{m-3}(x_s)] [f^{m-2}(x_s), x_s] ([x_s, x_{s+1}])^k$$ of length $m-2$, we again reach a contradiction.  Therefore, $q \ne m-2$.  This forces $q = m-1$ and so, $f^{m-1}(x_s) \le x_t < f^{m-2}(x_s)$.  Consequently, we have $f^{m-1}(x_s) = x_t$ and $z < f^{m-2}(x_s)$.  If $f^{m-3}(x_s) > z$, then by considering the cycle $$[x_s, x_{s+1}] [z, f(x_s)] [z : f^2(x_s)] \cdots [z : f^{m-4}(x_s)] [x_{s+1}, f^{m-3}(x_s)] [x_s, x_{s+1}]$$ of length $m-2$, we reach a contradiction.  Thus, $f^{m-1}(x_s) = x_t < x_{t+1} \le f^{m-3}(x_s) < z < f^{m-2}(x_s)$.    

Now we present the second argument.  It is clear that there is a {\it smallest} integer $1 \le r \le m-1$ such that both $f^r(x_s)$ and $f^{r+1}(x_s)$ lie on the same side of $z$.  If $r = 1$, then we consider the cycle $[x_s, x_{s+1}] [x_{s+1}, f(x_s)] ([x_s, x_{s+1}])^2$.  If $1 < r \le m-3$, then we consider the cycle $$[x_s, x_{s+1}] [z, f(x_s)] [x_s : f^2(x_s)] [f(x_s) : f^3(x_s)] \cdots [f^{r-2}(x_s) : f^r(x_s)] ([x_s, x_{s+1}])^{m-2-r}$$ of length $m-2$.  In either case, we have a periodic point of smaller odd period other than fixed point which is a contradiction.  So, $r \ge m-2$.  In particular, $(x_{t+1} \le ) \, f^i(x_s) < z < f^j(x_s)$ for all even $0 \le i \le m-3$ and all odd $1 \le j \le m-2$.  Since $q = m-2$ or $m-1$, this forces $r = m-1$ and so, we have $f^{m-1}(x_s) = x_t$ and $f^{m-1}(x_s) = x_t < x_{t+1} \le f^{m-3}(x_s) < z < f^{m-2}(x_s)$.  

Here is the third argument.  If $f(x_{\hat t}) \ge x_{s+1}$ for some $s+1 < \hat t \le m$, then since $q = m-2$ or $m-1$, it is clear that $x_{\hat t} \ne f^{m-2}(x_s), f^{m-1}(x_s)$.  So, $x_{\hat t} = f^r(x_s)$ for some $1 \le r \le m-3$.  If $r = 1$, we consider the cycle $[x_s, x_{s+1}]$ $[x_{s+1}, f(x_s)]$ $([x_s, x_{s+1}])^2$.  Otherwise, we consider the cycle $$[x_s, x_{s+1}] [z, f(x_s)] [z : f^2(x_s)] \cdots [z : f^{r-1}(x_s)] [x_{s+1}, f^r(x_s)] ([x_s, x_{s+1}])^{m-2-r}$$ of length $m-2$.  In either case, we have a periodic point of smaller odd period other than fixed point which is a contradiction.  Therefore, $f$ maps all $x_i$, $s+1 \le i \le m$, to the left side of $z$.  Since $f$ maps all $x_j$, $t+1 \le j \le s$ to the right side of $z$ and since $x_{t+1} \le f^i(x_s)$ for all $0 \le i \le m-3$, we obtain that $x_{t+1} \le f^{i_1}(x_s) < z < f^{i_2}(x_s)$ for all even $0 \le i_1 \le m-3$ and all odd $1 \le i_2 \le m-2$.  Thus $q = m-1$ and $f^{m-1}(x_s) = x_t$.  Consequently, $f^{m-1}(x_s) = x_t < x_{t+1} \le f^{m-3}(x_s) < z < f^{m-2}(x_s)$.  That $P$ is actually a \v Stefan cycle can be argued as that of Case 1 in section 8.

In the following, we present two proofs of (1) of Sharkovsky's theorem.  We also introduce a new general doubling operator and use it to construct new examples for (2) and (3).

\section{Proofs of Sharkovsky's theorem}
If $f$ has a period-$m$ point with $m \ge 3$ and odd, then it follows from (b) that $f$ has a period-$(m+2)$ point and, from (c) that $f$ has a period-$(2 \cdot 3)$ point.   

If $f$ has a period-$(2 \cdot m)$ point with $m \ge 3$ and odd, then, by Lemma 3(1), $f^2$ has a period-$m$ point.  It follows from the above  (or by (b) and (c)) that $f^2$ has a period-$(m+2)$ point and a period-$(2 \cdot 3)$ point.  If $f^2$ has a period-$(m+2)$ point, then, by Lemma 3(2), $$f \,\,\, \text{has either a period-}(m+2) \,\,\, \text{point or a period-}(2 \cdot (m+2)) \,\,\, \text{point}.$$  If $f$ has a period-$(m+2)$ point, then it follows from (c) that $f$ has a period-$(2 \cdot (m+2))$ point.  In either case, $f$ has a period-$(2 \cdot (m+2))$ point.  On the other hand, if $f^2$ has a period-$(2 \cdot 3)$ point, then, by Lemma 3(2), $f$ has a period-$(2^2 \cdot 3)$ point.  This shows that if $f$ has a period-$(2 \cdot m)$ point with $m \ge 3$ and odd, then $f$ has a period-$(2 \cdot (m+2))$ point and a period-$(2^2 \cdot 3)$ point.

Now if $f$ has a period-$(2^k \cdot m)$ point with $m \ge 3$ and odd and if $k \ge 2$, then, by Lemma 3(1), $f^{2^{k-1}}$ has a period-$(2 \cdot m)$ point.  It follows from the previous paragraph that $f^{2^{k-1}}$ has a period-$(2 \cdot (m+2))$ point and a period-$(2^2 \cdot 3)$ point.  So, by Lemma 3(2), $f$ has a period-$(2^k \cdot (m+2))$ point and a period-$(2^{k+1} \cdot 3)$ point.  

Furthermore, if $f$ has a period-$(2^i \cdot m)$ point with $m \ge 3$ and odd and if $i \ge 0$, then, by Lemma 3(1), $f^{2^i}$ has a period-$m$ point.  For each $\ell \ge i$, by Lemma 3(1), $f^{2^\ell} = (f^{2^i})^{2^{\ell-i}}$ has a period-$m$ point and so, by (c), $f^{2^\ell}$ has a period-6 point.  Thus, by Lemma 3(1), $f^{2^{\ell +1}}$ has a period-3 point and hence, by (a), has a period-2 point.  This implies, by Lemma 3(2), that $f$ has a period-$2^{\ell+2}$ point for each $\ell \ge i$.  

Finally, if $f$ has a period-$2^k$ point for some $k \ge 2$, then, by Lemma 3(1), $f^{2^{k-2}}$ has a period-4 point.  By (a), $f^{2^{k-2}}$ has a period-2 point.  By Lemma 3(2), $f$ has a period-$2^{k-1}$ point and hence, by induction, $f$ has a period-$2^j$ point for each $j = 1, 2, \cdots, k-2$.  Furthermore, it follows from (a) that $f$ has a fixed point.  This completes the proof of (1).

We now present a different proof of (1) which is inspired by the proof in {\bf\cite{bh}} and is based on the fact (proved in section 9) that if a period-$m$ orbit of $f$ contains a point that does not "switch sides" then $f$ has a period-$n$ point for each $n$ such that $n \prec m$.  We shall prove (1) by induction on $m$.  If $m =1$, there is nothing to prove. Suppose, for any continuous map from any nondegenerate compact interval into itself, if it has a period-$k$ point with $k < m$ then it also has a period-$\ell$ point for each $\ell \prec k$.  Now let $P$ be a period-$m$ orbit of $f$.  If $P$ contains a point that does not switch sides, then we are done.  Otherwise, assume that all points of $P$ switch sides and so $m$ is even.  Write $P = \{ x_1, x_2, \cdots, x_m \}$ with $x_1 < x_2 < \cdots < x_m$.  Then we have $f(\{x_1, x_2, \cdots, x_{m/2} \}) = \{ x_{(m/2)+1}, x_{(m/2)+2}, \cdots, x_m \}$ and $f(\{ x_{(m/2)+1}, x_{(m/2)+2}, \cdots, x_m \}) = \{ x_1, x_2, \cdots, x_{m/2} \}$.  Define a map $g(x)$ in $C^0(I, I)$ by putting $g(x) = \max \{ x_{(m/2)+1}, f(x) \}$ for $\min I \le x \le x_{m/2}$, $g(x) = \min \{ x_{m/2}, f(x) \}$ for $x_{(m/2)+1} \le x \le \max I$, and $g(x) = f(x)$ for $x_{m/2} \le x \le x_{(m/2)+1}$.  Then $g([\min I, x_{m/2}]) \subset [x_{(m/2)+1}, \max I]$ and $g([x_{(m/2)+1}, \max I]) \subset [\min I, x_{m/2}]$.  It is clear that $P \cap [\min I, x_{m/2}]$ is a period-$(\frac m2)$ orbit of $g^2$ in $[\min I, x_{m/2}]$ and $\frac m2 < m$.  By induction hypothesis, $g^2$ has a period-$i$ orbit $Q_i$ in $[\min I, x_{m/2}]$ for each $1 \le i$ and $i \prec \frac m2$.  Consequently, $Q_i \cup g(Q_i)$ is a period-$(2i)$ orbit of $g$ and also of $f$ and $2i \prec m$.  Since it is clear that $f$ has a fixed point, this completes the induction and hence the proof of (1).  On the other hand, we can also consider the continuous map $h(x)$ defined by putting $h(x) = \max \{ z, f(x) \}$ for $\min I \le x \le z$; and $h(x) = \min \{ z, f(x) \}$ for $z \le x \le \max I$.  Then $h^2$ is a continuous map from $[\min I, z]$ into itself and has a periodic point of least period $\frac m2 \, (< m)$ in $[\min I, z]$.  By induction hypothesis, $h^2(x)$ has a period-$j$ point in $[\min I, z]$ for each $1 < j$ and $j \prec \frac m2$.  So, discussed as above, $h$ has a period-$(2j)$ point (which turns out to be a period-$(2j)$ point of $f$ too) for each $1 < j$ and $2j \prec m$.  As for the existence of period-2 points for $f$, we note that, since $h^2$ has a period-2 orbit, say $\{ c, d \}$ with $c < d$ in $[\min I, z]$, $h^2$ has a fixed point in $(c, d) \, (\subset [\min I, z))$ which turns out to be a period-2 point of $h$ and hence of $f$.  This proves (1).  The following is an argument without referring to the map $g(x)$ or the map $h(x)$: If $m \ge 3$ is odd, then we are done.  If $m \ge 4$ is even, then $f^2$ has a period-$(\frac m2)$ point.  Since $\frac m2 < m$, it follows from the induction hypothesis that $f^2$ has a period-$k$ point for each $1 < k$ and $k \prec \frac m2$.  If $k$ is even, then $f$ has a period-$(2k)$ point.  If $k$ is odd, then, either $f$ has a period-$(2k)$ point or $f$ has a period-$k$ point which in turn implies that $f$ has a period-$(2k)$ point.  Therefore, in the case of even $m$, no matter whether $k$ is odd or even, as long as $1 < k$ and $2k \prec m$, $f$ has a period-$(2k)$ point.  In particular, $f$ has a period-4 point which, in turn, implies that $f$ has a period-2 point.  Finally, it is clear that $f$ has a fixed point.  This completes the induction and hence the proof of (1).  

As for the proofs of (2) and (3), it suffices (cf. {\bf\cite{bh}}) to consider the tent map $T(x) = 1 - |2x - 1|$ and the truncated tent map $\hat T_{a,b}(x)$, where $0 < a < b < 1$, defined on $[0, 1]$ by 
$$
\hat T_{a,b}(x) = \begin{cases}
               b, & \text{if} \,\,\, T(x) > b; \cr
               T(x), &  \text{if} \,\,\, a \le T(x) \le b; \cr
               a, & \text{if} \,\,\, T(x) < a.\cr
               \end{cases}
$$
The relationship between $T(x)$ and $\hat T_{a,b}(x)$ is that the periodic orbits of $\hat T_{a,b}(x)$ are also periodic orbits of $T(x)$ with the same periods and, conversely, the periodic orbits of $T(x)$ which lie entirely in the interval $[a, b]$ are also periodic orbits of $\hat T_{a,b}(x)$ with the same periods.  Consequently, if $Q_k$ is a period-$k$ orbit of $T(x)$, then it is also a period-$k$ orbit of $\hat T_{\min Q_k, \max Q_k}(x)$.  By (1), $\hat T_{\min Q_k, \max Q_k}(x)$ has a period-$\ell$ orbit for each $\ell$ with $\ell \prec k$.  In other words, the interval $[\min Q_k, \max Q_k]$ contains a period-$\ell$ orbit of $T(x)$ for each $\ell$ with $\ell \prec k$.  Since, for each integer $k \ge 1$, the equation $T^k(x) = x$ has exactly $2^k$ distinct solutions in $[0, 1]$, $T(x)$ has {\it finitely many} period-$k$ orbits.  Among these {\it finitely many} period-$k$ orbits, let $$P_k \,\,\, \text{be one with the {\it smallest} diameter} \,\,\, \max P_k - \min P_k.$$  For each $x$ in $I$, let $\hat T_k(x) = \hat T_{a_k,b_k}(x)$, where $a_k = \min P_k$ and $b_k = \max P_k(x)$.  Then it is easy to see that, for each $k \ge 1$, $\hat T_k(x)$ has exactly one period-$k$ orbit (i.e., $P_k$) but has no period-$j$ orbit for any $j$ with     $k \prec j$ in the Sharkovsky ordering.  This confirms (2).  

Clearly, $T(x)$ has a unique period-2 orbit, i.e., $\{ \frac 25, \frac 45 \}$.  For every periodic orbit $P$ of $T(x)$ with least period $\ge 3$, it follows from (a) that $\hat T_{\min P, \max P}(x)$ has a period-2 orbit.  So, $\min P \le \frac 25 < \frac 45 \le \max P$.  Now let $Q_3$ be any period-3 orbit of $T(x)$ of smallest diameter.  Then $[\min Q_3, \max Q_3]$ contains finitely many period-6 orbits of $T(x)$.  If $Q_6$ is one of smallest diameter, then $[\min Q_6, \max Q_6]$ contains finitely many period-12 orbits of $T(x)$.  We choose one, say $Q_{12}$, of smallest diameter and continue the process inductively.  Let $$q_0 = \sup\{\min Q_{2^n \cdot 3} : n \ge 0\} \,\,\, \text{and} \,\,\, q_1 = \inf\{\max Q_{2^n \cdot 3} : n \ge 0\}$$and let $\hat T_\infty(x) = \hat T_{q_0,q_1}(x)$ for all $0 \le x \le 1$.  If $\hat T_\infty(x)$ had a period-$(2^n \cdot m)$ orbit for some $n \ge 0$ and some odd $m \ge 3$, then, by (1), $\hat T_\infty(x)$ has a period-$(2^{n+1} \cdot 3)$ orbit, say $\hat Q_{2^{n+1} \cdot 3}$.  Since $\hat Q_{2^{n+1} \cdot 3} \subset [q_0, q_1] \subsetneq [\min Q_{2^{n+1} \cdot 3}, \max Q_{2^{n+1} \cdot 3}]$, $\hat Q_{2^{n+1} \cdot 3}$ is also a period-$(2^{n+1} \cdot 3)$ orbit of $T(x)$ with {\it smaller} diameter than that of $Q_{2^{n+1} \cdot 3}$.  This is a contradiction.  So, $\hat T_\infty(x)$ has no periodic orbit of period not a power of 2.  On the other hand, for each $k \ge 0$, the map $T(x)$ has finitely many period-$2^k$ orbits.  If each such orbit had an {\it exceptional} point which is not in the interval $[q_0, q_1]$, then it is clear that we can find an $n \ge 1$ such that the interval $[\min Q_{2^n \cdot 3}, \max Q_{2^n \cdot 3}]$ contains none of these {\it exceptional} points which implies that $[\min Q_{2^n \cdot 3}, \max Q_{2^n \cdot 3}]$ contains no period-$2^k$ {\it orbits} of $T(x)$.  Consequently, the map $\hat T_{s_n,t_n}(x)$, where $s_n = \min Q_{2^n \cdot 3}, t_n = \max Q_{2^n \cdot 3}$, has no period-$2^k$ {\it orbits} and yet it has a period-$(2^n \cdot 3)$ orbit, i.e., $Q_{2^n \cdot 3}$. This contradicts (1).  Therefore,the map $\hat T_\infty(x)$ is an example for (3).

On the other hand, there is a classical way of constructing examples for (2) and (3) through the doubling operator $\hat H_a(h)$ {\bf{\cite{al, st}}} (see Figure 1) which is a continuous map from $[0, 1]$ into itself defined, for any fixed $0 < a < 1/2$ and any continuous map $h$ from $[0, 1]$ into itself, by
$$
(\hat H_a(h))(x) =  \begin{cases}
               ah(x/a) + (1-a), & 0 \le x \le a, \cr
               \text{decreasing on} \,\, [a, 1-a],  \cr
               x - (1-a), & 1-a \le x \le 1. \cr
       \end{cases}
$$

\begin{figure}[htb]
\centerline{\epsfig{file=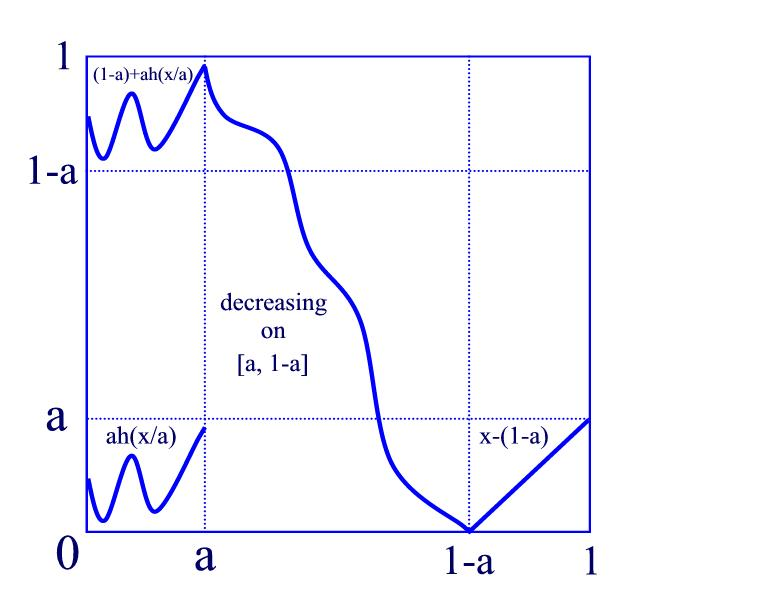,width=8cm,height=6cm}}
\caption{The graph of the classical doubling operator $\hat H_a(h)$ of the map $h_{|[0, 1]}$ on $[0, 1]$.}
\end{figure}

In the next section, we shall introduce a general doubling operator which includes the above classical doubling operator as a special case.

\section{A general doubling operator}
In this section, we introduce a general doubling operator which includes the above classical doubling operator as a special case.  Let $a$ and $b$ be two real numbers such that $0 < a < b < 1$.  Let $\psi(x)$ be a homeomorphism from $[0, a]$ onto itself and let $\phi(x)$ be a homeomorphism from $[0, a]$ onto $[b, 1]$.  For any continuous map $f$ from $[0, a]$ into itself, we let the doubling operator $(\Psi_a(f))(x)$ be the continuous map from $[0, 1]$ into itself defined by
$$
(\Psi_a(f))(x) =  \begin{cases}
               \phi(f(\psi(x))), & 0 \le x \le a, \cr
               \text{decreasing on} \,\, [a, b],  \cr
               (\psi^{-1} \circ {\phi}^{-1})(x), & b \le x \le 1. \cr
       \end{cases}
$$

It is clear that $(\Psi_a(f))([0, a]) \subset [b, 1]$ and $(\Psi_a(f))([b, 1])$ $\subset [0, a]$ and $(\Psi_a(f))^{2n}(x) = (\psi^{-1} \circ f^n \circ \psi)(x)$ on $[0, a]$.  So, $x_0 \in [0, a]$ is a period-$(2n)$ point of $\Psi_a(f)$ if and only if $\psi(x_0)$ is a period-$n$ point of $f$.  Furthermore, since $\Psi_a(f)$ is decreasing on $[a, b]$, $\Psi_a(f)$ can only have a fixed point, and maybe some period-$2$ points in $[a, b]$.  However, $\Psi_a(f)$ already has period-$2$ points in $[0, a]$, i.e., the fixed points of $f$.  Therefore, we obtain that $\{ m : \Psi_a(f)$ has a period-$m$ point in $[0, 1]\} = \{ m : \Psi_a(f)$ has a period-$m$ point in $[0, a]\} \cup \{ 1 \} = \{ 2n : f$ has a period-$n$ point in $[0, a] \}$ $\cup \{ 1 \}$.  Let $L_a(x) = ax, x \ge 0$.  For any continuous map $g$ from $[0, 1]$ into itself, the topologically conjugate map $(L_a \circ g \circ (L_a)^{-1})(x) = ag(x/a)$ is a continuous map from $[0, a]$ into itself.  So, we can let $(\hat \Psi_a(g))(x) = (\Psi_a(L_a \circ g \circ (L_a)^{-1}))(x), 0 \le x \le 1$.  Then $\hat \Psi_a$ doubles the periods of all periodic points of any continuous map from $[0, 1]$ into itself.  By taking $0 < a < 1/2 < b = 1-a < 1$, $\psi(x) = x$ and $\phi(x) = x + 1-a$, we obtain the classical doubling operator $\hat H_a$ as introduced in the previous section.  In the following, we take three more choices of $\phi$ and $\psi$, for $0 < a < 1/2 < b = 1-a < 1$, to obtain three more concrete doubling operators whose actions can be easily described in geometric ways. 

If, for $0 \le x \le a$, we let $\psi(x) = x$ and $\phi(x) = 1 - x$, then we have the following doubling operator $G_a(f)$ (see Figure 2) which is defined, for any continuous map $f$ from $[0, a]$ into itself, by
$$
(G_a(f))(x) =  \begin{cases}
               1 - f(x), & 0 \le x \le a, \cr
               \text{decreasing on} \,\, [a, 1-a], \cr
               1 - x, & 1-a \le x \le 1. \cr
       \end{cases}
$$

\begin{figure}[htb]
\centerline{\epsfig{file=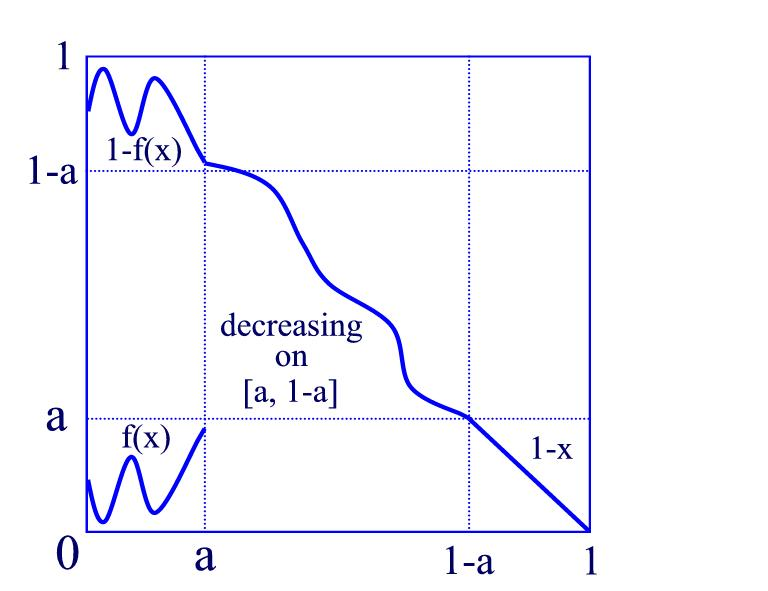,width=8cm,height=6cm}}
\caption{The graph of the doubling operator $G_a(f)$ of the map $f_{|[0, a]}$ on $[0, 1]$.}
\end{figure}
Note that the action of $G_a(f)$ on $[0, a]$ is reflecting the graph of $y = f(x)$ symmetrically with respect to the horizontal line $y = 1/2$ and, on $[1-a, 1]$ is reflecting the diagonal $y = x$ symmetrically with respect to the horizontal line $y = 1/2$.  For any continuous map $g$ from $[0, 1]$ into itself, the topologically conjugate map $ag(x/a)$ is a continuous map from $[0, a]$ into itself, let $(\hat G_a(g))(x) = G_a(L_a \circ g \circ (L_a)^{-1})(x)$, where $L_a(x) = ax, x \ge 0$.  Then $\hat G_a$ doubles the periods of all periodic points of any continuous map from $[0, 1]$ into itself.

If, for $0 \le x \le a$, we let $\psi(x) = a - x$ and $\phi(x) = 1 - a + x$, then we have the following doubling operator $F_a(f)$ (see Figure 3) which is defined, for any continuous map $f$ from $[0, a]$ into itself, by
$$
(F_a(f))(x) =  \begin{cases}
               1-a + f(a-x), & 0 \le x \le a, \cr
               \text{decreasing on} \,\, [a, 1-a], \cr
               1 - x, & 1-a \le x \le 1. \cr
       \end{cases}
$$
Note that the action of $F_a(f)$ on $[0, a]$ is reflecting the graph of $y = f(x)$ symmetrically with respect to the vertical line $x = a/2$ and then pushing the resulting graph up $1-a$ units and, on $[1-a, 1]$ is reflecting the diagonal $y = x$ symmetrically with respect to the vertical line $x = 1 - a/2$ and then pulling the resulting graph down $1-a$ units.  For any continuous map $g(x)$ from $[0, 1]$ into itself, the doubling operator $(\hat F_a(g))(x)$ is defined similarly.  

\begin{figure}[htb]
\centerline{\epsfig{file=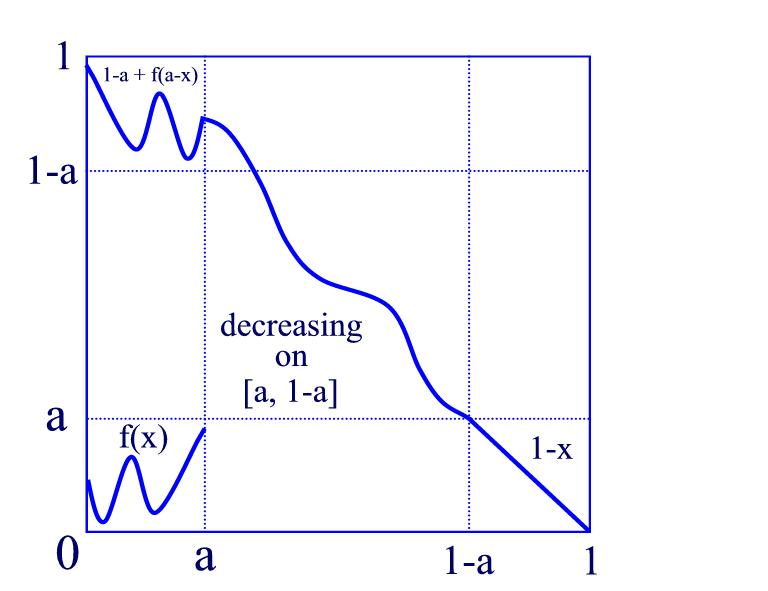,width=8cm,height=6cm}}
\caption{The graph of the doubling operator $F_a(f)$ of the map $f_{|[0, a]}$ on $[0, 1]$.}
\end{figure}

If, for $0 \le x \le a$, we let $\psi(x) = a - x$ and $\phi(x) = 1 - x$, then we have the following doubling operator $E_a(f)$ (see Figure 4) which is defined, for any continuous map $f$ from $[0, a]$ into itself, by
$$
(E_a(f))(x) =  \begin{cases}
               1 - f(a-x), & 0 \le x \le a, \cr
               \text{decreasing on} \,\, [a, 1-a], \cr
               x - (1-a), & 1-a \le x \le 1. \cr
       \end{cases}
$$
Note that the action of $E_a(f)$ on $[0, a]$ is reflecting the graph of $y = f(x)$ symmetrically with respect to the point $(a/2, 1/2)$ and, on $[1-a, 1]$ is reflecting the diagonal $y = x$ symmetrically with respect to the point $(1 - a/2, 1/2)$.  For any continuous map $g(x)$ from $[0, 1]$ into itself, the doubling operator $(\hat E_a(g))(x)$ is defined similarly.  

\begin{figure}[htb]
\centerline{\epsfig{file=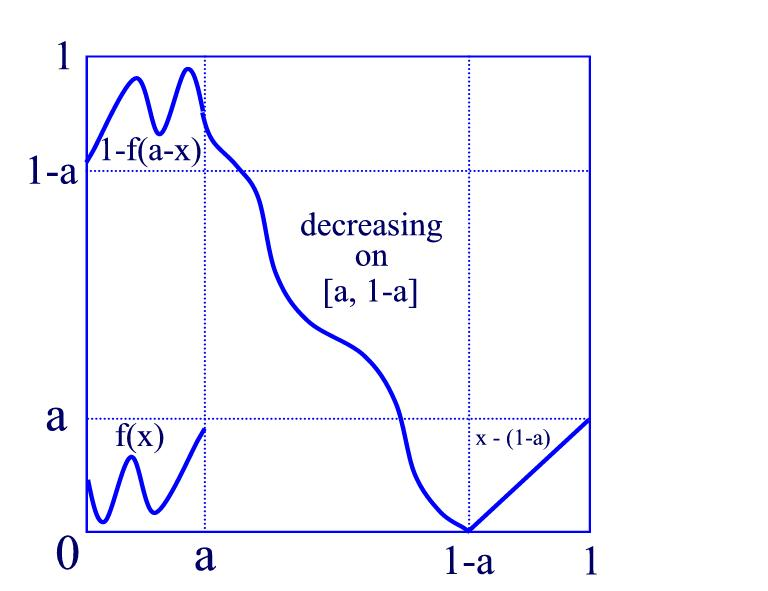,width=8cm,height=6cm}}
\caption{The graph of the doubling operator $E_a(f)$ of the map $f_{|[0, a]}$ on $[0, 1]$.}
\end{figure}

For the sake of completeness, we include the classical \v Stefan map $f_n(x)$ which has a perodic orbit of least period $2n+1$, but no periodic orbits of any smaller odd periods other than fixed points.  For every integer $n \ge 2$, let $x_i, 1 \le i \le 2n+1$, be $2n+1$ distinct points in $[0, 1]$ such that $0 = x_1 < x_2 < \cdots < x_{2n} < x_{2n+1} = 1$.  Let $f_n(x)$ be the continuous map from $[x_1, x_{2n+1}] \, (= [0, 1])$ onto itself defined by putting $f_n(x_1) = x_{n+1}$; $f_n(x_i) = x_{2n+3-i}$ for $2 \le i \le n+1$; $f_n(x_j) = x_{2n+2-j}$ for $n+2 \le j \le 2n+1$; and by linearity on each interval $[x_k, x_{k+1}]$, $1 \le k \le 2n$.  Then it is easy to see that the set $\{ x_i : 1 \le i \le 2n+1 \}$ forms a \v Stefan cycle of $f_n$ with least period $2n+1$ and $f_n^{2n-1}([x_1, x_2]) = [x_2, x_{2n+1}]$.  So, $f_n^{2n-1}([x_1, x_2]) \cap [x_1, x_2) = \emptyset$.  That is, $[x_1, x_2]$ contains no period-$(2n-1)$ point of $f_n$.  Furthermore, since $f_n$ is strictly decreasing on $[x_2, x_{2n+1}]$, any periodic orbit of least period $> 2$ must have at least one point in $[x_1, x_2]$ and hence $f_n$ cannot have period-$(2n-1)$ points.  Therefore, $f_n$ has a  period-$(2n+1)$ orbit but no period-$(2n-1)$ orbits.  On the other hand, on $[0, 1]$, let $g(x) \equiv 0$ be the constant map and let $h(x) = 1/2 + 2x$ for $0 \le x \le 1/4$; $h(x) = -2x + 3/2$ for $1/4 \le x \le 1/2$; and $h(x) = 1 - x$ for $1/2 \le x \le 1$.  Then $g$ has fixed points but no period-2 points and $h$ has period-6 points (for example, the point $x = 1/2 - 1/7$) but no periodic points of any odd periods $> 1$.  By applying successively any one of the above doubling operators $\hat \Psi_a$, $\hat H_a$, $\hat G_a$, $\hat F_a$ or $\hat E_a$ to $f_n$, $g$, and $h$ respectively, we obtain examples for (2).  

Finally, we present an example for (3) which is different from the classical one as described in {\bf\cite{al, st}}.  Let $< a_i >_{i \ge 1}$, $< b_i >_{i \ge 1}$ be any two infinite sequences of numbers in $(0, 1/2)$ and let $< \al_i >_{i \ge 1}$ be any infinite sequence of $0$'s and $1$'s.  For $i \ge 1$, let 
$$
\Phi_{\al_i} = \begin{cases}
               \hat G_{a_i}, & \text{if} \,\, \al_i = 0, \\
               \hat H_{b_i}, & \text{if} \,\, \al_i = 1. \\
               \end{cases}
               \qquad \text{and} \qquad
c_i = \begin{cases}
               a_i, & \text{if} \,\, \al_i = 0, \\
               b_i, & \text{if} \,\, \al_i = 1. \\
               \end{cases}
$$
Then, for any fixed $n > 2$ and any continuous map $\phi$ from $[0, 1]$ into itself, since we have $(\Phi_{\al_n}(\Phi_{\al_{n+1}}(\phi)))(x) = (\Phi_{\al_n}(\phi))(x)$ on $[1-c_n, 1]$, the following hold (may not hold if in the definition of $\Phi_{\al_i}$, $\hat G_{a_i}$ or $\hat H_{b_i}$ is replaced by $\hat E_{a_i}$, $\hat F_{a_i}$, $\hat E_{b_i}$, or $\hat F_{b_i}$), $$(\Phi_{\al_{n-1}}(\Phi_{\al_n}(\Phi_{\al_{n+1}}(\phi))))(x) = (\Phi_{\al_{n-1}}(\Phi_{\al_n}(\phi)))(x) \quad  \text{on} \quad [c_{n-1}(1-c_n), 1]$$ and $$|(\Phi_{\al_{n-1}}(\Phi_{\al_n}(\Phi_{\al_{n+1}}(\phi))))(x) - (\Phi_{\al_{n-1}}(\Phi_{\al_n}(\phi)))(x)| < c_{n-1} \quad \text{on} \quad [0, c_{n-1}(1-c_n)].$$  By induction, we obtain that, on $[c_1c_2c_3 \cdots c_{n-2}c_{n-1}(1-c_n), 1]$, 
$$
(\Phi_{\al_1}(\Phi_{\al_2}(\cdots(\Phi_{\al_n}(\Phi_{\al_{n+1}}(\phi))\cdots)(x) = (\Phi_{\al_1}(\Phi_{\al_2}(\cdots(\Phi_{\al_n}(\phi))\cdots)(x), 
$$
and, on $[0, c_1c_2c_3 \cdots c_{n-2}c_{n-1}(1-c_n)]$,
$$
|(\Phi_{\al_1}(\Phi_{\al_2}(\cdots(\Phi_{\al_n}(\Phi_{\al_{n+1}}(\phi))\cdots)(x) - (\Phi_{\al_1}(\Phi_{\al_2}(\cdots(\Phi_{\al_n}(\phi))\cdots)(x)| < \Pi_{i=1}^{n-1} c_i < 1/{2^{n-1}}.
$$
Thus, the sequence $\Phi_{\al_1}(\phi)$, $\Phi_{\al_1}(\Phi_{\al_2}(\phi))$, $\Phi_{\al_1}(\Phi_{\al_2}(\Phi_{\al_3}(\phi)))$, $\cdots$, which double the periods of all periodic points of $\phi(x)$ successively, converges uniformly to a continuous map $\Phi_{\al}$, where $\alpha = \alpha_1\alpha_2\cdots$, on $[0, 1]$ (which is {\it independent} of $\phi$).  It is easy to see that $\Phi_{\al}$ is an example for (3).  Since there are uncountably many $\al$'s, we have uncountably many examples $\Phi_\al$ for (3).  Figure 5 is such an example with (i) $a_i = b_i = 1/3$, $i \ge 1$; (ii) $\beta = \beta_1\beta_2\beta_3 \cdots = \overline {01} = 010101 \cdots$; and (iii) both $\hat G_{a_i}$ and $\hat H_{b_i}$, $i \ge 1$, are linear on $[1/3, 2/3]$. 

\begin{figure}[htb]
\centerline{\epsfig{file=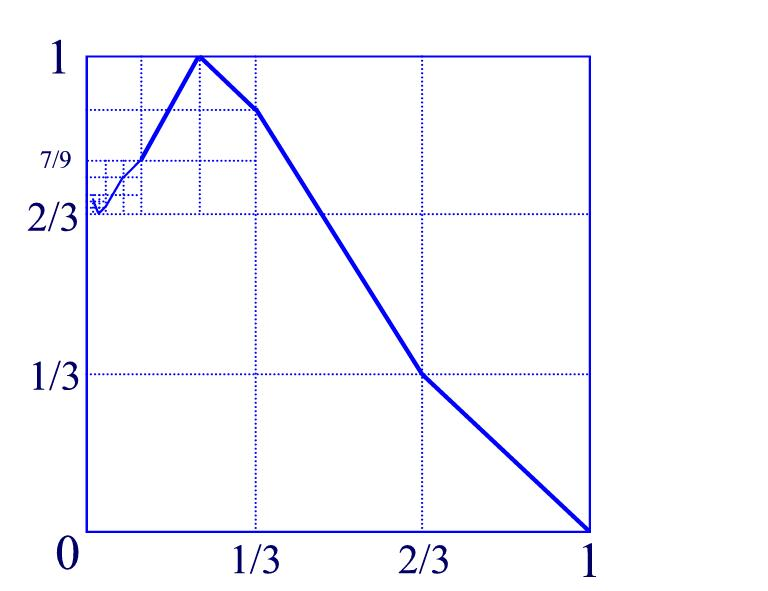,width=8cm,height=6cm}}
\caption{The graph of the map $\Phi_{\overline{01}}$. Note that $\Phi_{\overline{01}}(0) = 7/10$.}
\end{figure}

\noindent
{\bf Remark.}
In the above construction of $\Phi_\al$ with $a_i = b_i = 1/3$, $i \ge 1$, if $\al = \overline{1} = 111 \cdots$, then we have the classical example as introduced in {\bf{\cite{al}}} with $\Phi_{\overline 1}(0) = 1$.  On the other hand, for any $\al \ne \overline 1$, there is a $k \ge 1$ such that $\al_k = 0$.  Then because of the turning upside down property of $\hat G_{a_k} (= \Phi_{\al_k})$, we obtain that $0 < \Phi_{\al}(0) < 1$.  For example, for $\beta = \overline{01}$ considered above, we have $\Phi_\beta(0) = 7/10$.  Similarly, if $a_i \equiv c$, $i \ge 1$, where $c$ is a constant in $(0, 1/2)$, and if $\gamma = \overline{0}$, then we have $\Phi_\gamma(0) = 1/(1+c)$ $\in (0, 1)$.

\end{document}